%% file: self-consisbnds-main.tex
\documentclass{article}
\usepackage{amssymb}
\usepackage{amsmath}
\usepackage{graphicx}
\usepackage{xspace}

\newtheorem{theorem}{Theorem}
\newtheorem{definition}{Definition}
\newtheorem{lemma}[theorem]{Lemma}

\newtheorem{remark}[theorem]{Remark}

\newcommand{\dom }{{\rm dom}\,}

\newcommand{\gss}{GSS\xspace}
\newcommand{\VL}{{\rm\textbf{VC}}\xspace}
\newcommand{\NC}[1]{\textbf{N#1}}
\newcommand{\BD}[3]{#1_{\langle #2 \rangle \langle #3 \rangle}}

\def\qed{{\hfill{\vrule height5pt width3pt depth0pt}\medskip}}
\def\comment#1{{}}

\begin{document}
\begin{center}
{\Large \bf  Self-consistent bounds method for dissipative PDEs}

 \vskip 0.5cm
{\large Daniel Wilczak} and {\large Piotr Zgliczy\'nski}\footnote{Work of D.W. and P.Z. was supported by National Science Center (NCN) of Poland under project No. UMO-2016/22/A/ST1/00077} 
 \vskip 0.2cm
  Jagiellonian University, Faculty of Mathematics and Computer Science, \\
  \L ojasiewicza 6, 30--348  Krak\'ow, Poland \\ \texttt{e-mail:\,\{Daniel.Wilczak,Piotr.Zgliczynski\}@uj.edu.pl}
\vskip 0.5cm

\today
\end{center}

\begin{abstract}
 We discuss the method of self-consistent bounds for dissipative PDEs with periodic boundary conditions. We prove convergence theorems for a class of dissipative PDEs, which constitute a theoretical basis of a general framework for construction of an algorithm that computes bounds for the solutions of the underlying PDE and its dependence on initial conditions. 
 
 We also show, that the classical examples of parabolic PDEs including Kuramoto-Sivashinsky equation and the Navier-Stokes on the torus fit into this framework. 
\end{abstract}

\input introscb.tex

\input gcp-norms.tex

\input lin-maps-gss.tex

\input lineqestm.tex

\input c0conver.tex

\input c1conver.tex

\input isolation.tex

\input convercond.tex

\input varconvcond.tex

\input ref.tex
\end{document}

%% file: introscb.tex
\section{Introduction}

The main motivation for this research is to provide a theoretical framework for construction of $\mathcal C^1$-like algorithm for a class of dissipative PDEs. By a $\mathcal C^1$ algorithm we understand computation that provides guaranteed bounds on forward trajectories of underlying semi-flow and bound on the solutions to the associated variational equations.

The method of self-consistent bounds introduced in \cite{ZM}  lead to design  of rigorous integrators for a class of dissipative PDEs  \cite{C,WZ,ZKSper,ZKS3}, which  proved to be very efficient in studying dynamics of evolutionary dissipative PDEs \cite{BKZ,WZ,CW,CZ,Z,Z1,ZKSper,ZKS3}.

The aim of this paper is threefold. First, we revise the method of self-consistent bounds and reformulate it in an abstract setting. The statements are proved under very general assumptions about norms used, sets on which the existence of local semi-flow is proved and, the most important, assumptions on the vector field that guarantee uniform convergence of the solutions of finite-dimensional approximations (Galerkin projections) to the solutions of the original infinite-dimensional system. We will later refer to these assumptions by the letters \textbf{N} for norms, \textbf{S} for sets and \textbf{C} for convergence conditions.

The second aim is to formulate and prove results about convergence of the solutions to associated variational equations for Galerkin projections to the solutions of (properly understood) variational equations of infinite-dimensional systems. Here we will need additional convergence conditions for variational equations, we refer to as \textbf{VC}.

Finally we will show, that convergence conditions \textbf{C} and \textbf{VC} are not very restrictive, and they are naturally satisfied in a class of vector fields including widely studied models such as Kuramoto-Sivashinsky  \cite{KT,S}, Burgers \cite{B}, Brusselator (reaction-diffusion system) and Navier-Stokes \cite{T} equations. The crucial fact which make the proposed framework applicable is so-called \emph{isolation property} of the vector field, which informally means that the vector field is pointing inwards for all high modes leading to dissipative character of the system.

To illustrate the isolation property we consider a class of  PDEs of the following form
\begin{equation}
  u_t = L u + N\left(u,Du,\dots,D^ru\right), \label{eq:genpde-intro}
\end{equation}
where $u \in \mathbb{R}^n$,  $x \in \mathbb{T}^d=\left(\mathbb{R}\mod 2\pi\right)^d$, $L$ is a linear operator, $N$ is a
polynomial and by $D^s u$ we denote $s^{\text{th}}$ order derivative of
$u$, i.e. the collection of all spatial (i.e. with respect to variable $x$) partial derivatives of $u$ of
order $s$. 

We require, that the operator $L$ is diagonal in the Fourier basis
$\{e^{ikx}\}_{k \in \mathbb{Z}^d}$,
\begin{equation*}
  L e^{ikx}= -\lambda_k e^{ikx},
\end{equation*}
with
\begin{eqnarray*}
 L_* |k|^p &\leq& \lambda_k \leq  L^* |k|^p, \qquad \text{for all $|k| > K$ and  $K,L_*,L^* \geq 0$}, \\
    p &>& r.
\end{eqnarray*}
The assumption $p>r$ makes the equation dissipative,  we have a smoothing effect.


If the solutions are sufficiently smooth, the problem (\ref{eq:genpde-intro}) can be written as an infinite ladder of ordinary differential equations for the Fourier
coefficients in $u(t,x)=\sum_{k \in \mathbb{Z}^d} u_k(t) e^{i kx}$, as follows
\begin{equation}
  \frac{d u_k}{dt} = F_k(u)=-\lambda_k u_k + N_k\left(\{u_j\}_{j \in \mathbb{Z}^d}\right), \qquad \mbox{for all
  $k \in \mathbb{Z}^d$}. \label{eq:fueq}
\end{equation}

The crucial fact, which makes our approach  possible is the \emph{isolation property}, which will be defined in Section~\ref{sec:isolation}.
Here, informally we express it as follows.

\emph{Let}
\begin{equation*}
   W=\left\{ \{u_k\}_{k \in \mathbb{Z}^d}\,:\,  |u_k| \leq \frac{C}{|k|^s q^{|k|}}\right\}
\end{equation*}
\emph{where $q\geq 1$, $C>0$, $s\geq 0$, satisfy some inequalities \\
Then under some assumptions there  exists $K>0$, such that for $|k| > K$ there holds}
\begin{equation*}
 \mbox{if} \quad u \in W, \ |u_k|=\frac{C}{|k|^s q^{|k|}}, \qquad \mbox{then} \qquad  u_k \cdot F_k(u) <0.
\end{equation*}

When $q=1$ and $s\geq p+d+1$ this is established in   Lemma~\ref{lem:iso-mainVar} in Section~\ref{subsec:iso-main-var} (see also Theorem 3.1 in \cite{ZKS3}). In  \cite{WZ} we used sets of the type described above with $s=0$ and $q>1$. 

Projection of the set $W$ onto $k^{\mathrm{th}}$ mode is a closed disc (or interval) centred at zero and of radius
$r_k =\frac{C}{|k|^s q^{|k|}}$. Geometrically the isolation property means that if $|u_k|=r_k$ for some $k> K$ then the $k$-th component of the vector field is pointing inwards the set. As a consequence, the only way a trajectory may leave the set $W$ forward in time is by increasing some of leading modes $|u_k|$, $k\leq K$ above the threshold $r_k$. However, one has to be careful with what "inwards" means in the above statement because $W$ has empty interior as a compact subset in infinite dimensional space. On the other side it makes perfect sense for Galerkin projections.

This property is used   in  our approach to obtain a priori bounds for $u_k(h)$ for small $h>0$ and $|k|>K$, while the leading modes $u_k$ for $|k| \leq K$ are computed using  tools for rigorous integration of ODEs \cite{CAPDREVIEW,Lo,NJP} based on the interval arithmetics \cite{Mo}. Moreover, from the point of view of topological method, the isolation property is of crucial importance as it shows that on the tail we have the entry behaviour, which enable us to apply the finite dimensional tools from the dynamics like
the Conley index \cite{ZM}, the  fixed point index and  covering relations \cite{ZKS3}.

In the present work we focus on the issues related to the convergence of Galerkin projections of equations (\ref{eq:fueq}) and its variational equations.

The paper is organized as follows. In Section~\ref{sec:gss} we introduce the notion of good sequence space (\gss space) and  establish some of their properties.  \gss will be later the spaces in which our Fourier series belong. In Section~\ref{sec:estmlinEq} we recall the notion of the logarithmic norm of a matrix and prove a useful result regarding block decomposition for nonautonomous linear ODEs. This is a crucial result for efficient implementation of an algorithm based on the framework introduced in Section~\ref{sec:c0-conver} and Section~\ref{sec:c1-conver}. There we state and prove the main results about the convergence of solutions of finite-dimensional Galerkin projections to solutions of underlying PDE for the main system (Theorem~\ref{thm:limitLN}) and for variational system (Theorem~\ref{thm:c1conver}), respectively. In Section~\ref{sec:isolation} we formally define the isolation property and we prove that it implies the existence of a-priori bounds, which appear in the assumptions the convergence theorems. Finally, in Section~\ref{sec:dissipativePDEs} we show that the introduced framework applies to a class of dissipative DPEs, that is the assumptions on sets \textbf{S}, convergence conditions \textbf{C}, \textbf{VC} and the isolation property are satisfied for various choices of Banach spaces.

\subsection*{Notation}

By $\mathbb{Z}_+$ we denote the set of all positive integers.
Given two normed vector spaces $V,W$ by $\text{Lin}(V,W)$ we will denote the space of all bounded linear maps from $V$ to $W$.


%% file: gcp-norms.tex
\section{Sequence spaces with good norms}
\label{sec:gss}

The goal of this section is to define the sequence spaces on which we will work later in this paper. We will be dealing with sequence spaces $H \subset \mathbb{R}^{\mathbb{Z}_+}$ equipped with some norm $\|\cdot\|$. Put $e_i=(\delta_{ij})_{j\in\mathbb Z_+}$. By $e_j^* \in H^*$ we denote the dual form, that is $e^*_j(e_i)=\delta_{ij}$.

By $\pi_k$ we denote the projection onto $k$-th direction, i.e. $\pi_k\left(\sum w_j e_j\right)=w_ke_k$. For a nonempty set $J \subset \mathbb{Z}_+$ by $P_J$ we denote a projection defined by $P_J(w)=\sum_{i\in J}w_i e_i$.


For $n \in \mathbb{Z}_+$  by $H_n$ we denote a subspace spanned by
$\{e_j\}_{j \leq n}$. Put $P_n:=P_{\{j \leq n\}}$ and $Q_n=\mathrm{Id}-P_n$. By $\iota_n:H_n \to H$ we denote the embedding $H_n$ into $H$.

We define our standing assumptions \NC1--\NC5 on the space $(H,\|\cdot\|)$.

\begin{definition}
We will say that $H \subset \mathbb{R}^{\mathbb{Z}_+}$ with norm $\| \cdot \|$ is \emph{\gss} (good sequence space) if the following conditions are satisfied.
\begin{description}
\item[\NC 1] $H$ is a Banach space.
\item[\NC 2] $\forall w \in H\ w=\sum_i w_i e_i$.
\item[\NC 3] For all  $w \in H$ and for any $\alpha \in \{-1,1\}^{\mathbb{Z}_+}$ there holds $w^\alpha=\sum_i \alpha_i w_i e_i \in H$
and $\|w\|=\|w^\alpha\|$.
\item[\NC 4] $\|P_J w\| \leq \|w\|,  \quad \forall w \in H, \, \forall J \subset \mathbb{Z}_+$.
\item[\NC 5]  there exists  constant $G$, such that for all $w \in H$ and for all $i$ $|w_i| \leq G \|w\|$.
\end{description}
\end{definition}

\textbf{Examples:} $l_2$, $l_1$, $c_0$ (sequences converging to $0$) with the norm $\|\cdot\|_\infty$. However $l_\infty$ is not a \gss space, because \NC2 is not satisfied.
The space defined by convergence of $\sum_i |w_i|/2^i$ does not satisfy \NC 5.

\begin{lemma}
\label{lem:gcp-prop}
If $(H,\|\cdot\|)$  is \gss, then  for any $J \subset \mathbb{Z}_+$ the projection $P_J$ is continuous and for any $w \in H$ there holds
\begin{eqnarray}
    \|w\|&=&\sup_{n \in \mathbb{Z}_+} \|P_n w\|, \label{eq:nsupGalproj} \\
    \|(I-P_n) w\| &\to& 0, \quad \mbox{for $n \to \infty$}.  \label{eq:remto0}
\end{eqnarray}
\end{lemma}
\textbf{Proof:}
Continuity of $P_J$ follows immediately from \NC4. From \NC2 and \NC4 we have
$$\|w\|=\lim_{n\to\infty}\|P_n w^n\|\leq \sup_n\|P_nw\|\leq \|w\|,$$
which proves (\ref{eq:nsupGalproj}).
From \NC4 the sequence $\|P_nw\|$ is non-decreasing, hence $\lim_{n\to\infty}\|P_nw\|=\sup_{n}\|P_nw\|=\|w\|$. In consequence,
$$\|(I-P_n) w\|\leq \|w\|-\|P_nw\|\to 0.$$
\qed

\begin{lemma}
\label{lem:compt-gss}
  Assume that $H$ is \gss. Then $W \subset H$ is compact, iff it is closed, bounded and \emph{$W$ has uniform bounds on the tail}, i.e. for every $\epsilon >0$ there exists $N$, such that
  for all $n \geq N$ and for all $w \in W$ there holds $\|(I-P_n)w\| < \epsilon$.
\end{lemma}
\textbf{Proof:}
Implication $\Rightarrow$. Since $W$ is compact, it is also bounded and closed.  To prove the existence of uniform bound for the tail let us fix $\epsilon >0$
and define an open covering $\{U_n\}_{n \in \mathbb{Z}_+}$ of $W$ by setting $U_n=\{w \in H: \|(I-P_n)w\| < \epsilon\}$. It is indeed a covering due to (\ref{eq:remto0}) in Lemma~\ref{lem:gcp-prop}.
Observe that $U_n \subset U_{n+M}$ for any $n,M \in \mathbb{Z}_+$ (this a consequence of assumption \NC4).
From compactness of $W$ it follows that there exists $N$, such that $W \subset U_N \subset U_{N+M}$ for any $M \in \mathbb{Z}_+$.

Implication $\Leftarrow$.  Fix a sequence $c^k \in W$. We would like to prove that it has a convergent subsequence. From the boundedness of $W$ and the diagonal argument it follows that we can find a subsequence of $d^m=c^{k_m}$, such that for each $m \in \mathbb{Z}_+$ the sequence $(d^m_n)_{m\in\mathbb Z_+}$ is convergent to some $\bar{d}_n$. Therefore, for each $n$ the sequence $(P_nd^m)_{m\in\mathbb Z_+}$ satisfies the Cauchy condition. We will show that $d^m$ also satisfies the Cauchy condition.

Let us fix $\epsilon >0$. For $j\in{\mathbb Z_+}$ we have
\begin{eqnarray*}
  \|d^m - d^{m+j}\| \leq \|P_n d^m - P_n d^{m+j}\| + \|(I-P_n) d^m\| + \|(I-P_n) d^{m+j}\|.
\end{eqnarray*}
From the uniform bound on tail in $W$ we find $N$ such that $\|(I-P_N)w\| < \epsilon$ for $w \in W$, in particular for all $d^m$. Given fixed $N$ from the Cauchy condition for the sequence $(P_N d^m)_{m\in\mathbb Z_+}$ it follows that there exists $M$, such that for all $m \geq M$ and all $j \in \mathbb{Z}_+$ there holds $\|P_N d^m - P_N d^{m+j}\| < \epsilon$.

Hence we obtain for $m \geq M$ and all $j \in \mathbb{Z}_+$
 \begin{eqnarray*}
  \|d^m - d^{m+j}\| \leq 3 \epsilon.
\end{eqnarray*}
Therefore the sequence $d^m$ satisfies the Cauchy condition. Since $H$ is complete (assumption \NC1) there exists $d \in H$, such that $d^m \to d$. Since $W$
is closed, we conclude that $d \in W$.
\qed

From the proof of the above lemma it is  easy to infer the following result.
\begin{lemma}
\label{lem:convW-gss}
  Assume that $H$ is \gss  and $W \subset H$ is compact. Then, a sequence $c^k$ in $W$ is convergent to $c \in W$, iff for all $j\in \mathbb{Z}_+$ $\lim_{k \to \infty} c^k_j=c_j$.
\end{lemma}

\subsection{Linear forms on \gss}

The following lemma gives a characterization of continuous linear forms on \gss.
\begin{lemma}
\label{lem:formOnH}
Assume that $H$ is \gss and let $f:H \to \mathbb{R}$ be bounded linear functional. Put $f_j=f(e_j)$. Then, for every $a \in H$ there holds
\begin{eqnarray}
   \sum_{j \in \mathbb{Z}_+} |f_j| \cdot  |a_j| \leq \|f\| \cdot \|a\|  \label{eq:fkak-absconver}&\qquad \text{and}\\
   f(a)=\sum_{j \in \mathbb{Z}_+} f_j a_j.   \label{eq:funcLinonH}
\end{eqnarray}
\end{lemma}
\textbf{Proof:} Put $K=\|f\|$.
   Let $\bar{a} \in H$ be such that $|a_j|=|\bar{a}_j|$ and $f_j \bar{a}_j \geq 0$ for all $j$. From \NC3 we have $\bar{a} \in H$ and $\|a\|=\|\bar{a}\|$.  Since each term $f_j \bar{a}_j$ is nonnegative, using \NC4 we obtain
\begin{eqnarray*}
 K\|a\| = K \|\bar{a}\| \geq   K \|P_n \bar{a}\|  \geq    f(P_n \bar{a}) =   \sum_{j \leq n}  f_j  \cdot  \bar{a}_j =  \sum_{j \leq n}  |f_j|  \cdot  |a_j|
\end{eqnarray*}
hence we obtained (\ref{eq:fkak-absconver}).

Equation (\ref{eq:funcLinonH}) is a consequence of the continuity of $f$ and the assumption \NC2.

\qed

\subsection{Bounded linear self-maps on \gss}


\begin{theorem}
\label{thm:lin-gss}
Let $H$ be \gss space.
Assume that $V:H \to H$ is bounded linear map. Then there exists a collection of real numbers $\{V_{ij}\}_{i,j \in \mathbb{Z}_+^2}$, such that for each $a \in H$ there holds
\begin{eqnarray*}
  \sum_{j} |V_{ij}| \cdot |a_j| &\leq& \|\pi_i V\| \cdot \|a\|,  \quad i \in \mathbb{Z}_+ \\
  Va &=& \sum_{i} \left(\sum_{j}V_{ij}a_j\right)e_i.
\end{eqnarray*}
\end{theorem}
\textbf{Proof:}
We apply \NC2 and Lemma~\ref{lem:formOnH} to each $\pi_i V$, which are bounded linear forms.
\qed

In view of the above theorem we will often represent $V \in \mbox{Lin}(H,H)$ as a matrix  $V=\{V_{ij}\} \in \mathbb{R}^{\mathbb{Z}_+ \times \mathbb{Z}_+}$ and by $V_{\ast j}$ we will denote its $j$-th column.

%% file: lin-maps-gss.tex
\subsection{Galerkin projection of linear maps on \gss}

We will be interested in the following question: given linear maps $V^n: H_n \to H_n$ (for example the Jacobian matrices for the flow defined by Galerkin projections), does there exists a limit $V:H \to H$, under some assumptions, which  are reasonable in the context of dissipative PDEs.  To make sense of this question we treat map  $V^n:H_n \to H_n$  as a linear map from $H$ to $H$ by setting
$V^n(u)= \iota_n (V^n(P_n u)) \in H_n \subset H$.

Observe that in view of Theorem~\ref{thm:lin-gss} we can represent in basis $\{e_j\}$ each continuous linear map $V:H \to H$ by its matrix $\{V_{ij}\}_{i,j \in \mathbb{Z}_+}$. 
We also set $V_{*j}=Ve_j$, which can be seen as $j$-th column of matrix $V_{ij}$ representing $V$.

\begin{theorem}
\label{thm:Vn-weak-lim}
 Let $V^n:H \to H$, $n\in\mathbb Z_+$ be linear maps such that $\|V^n\| \leq K$ for some $K\geq 0$. Assume that
 \begin{itemize}
 \item for all $i,j \in \mathbb{Z}_+$ there exists $V_{ij}$, such that
 \begin{equation}
   \lim_{n\to \infty} V^n_{ij} = V_{ij} \label{eq:dfg-Vij-conv0}
 \end{equation}
 \item there exists a family of compact sets $W^j\subset H$, $j\in\mathbb Z_+$ such that for all $n$ there holds
 \begin{equation}
 V^n_{\ast j} \in W^j. \label{eq:Vnj-approribnd}
 \end{equation}
\end{itemize}

Then there exists a bounded linear map $V:H\to H$, such that
\begin{enumerate}
\item  for all $a \in H$ there holds $\lim_{n \to \infty} V^n(a)=Va$,
\item $\|V\| \leq K$,
\item  $\pi_iV(e_j)= V_{ij} e_i, \quad \forall i,j\in \mathbb{Z}_+$ and
\item for each $i$ and all $a \in H$ there holds $\sum_{j=1}^\infty |V_{ij}| \cdot |a_j| < \infty$.
\end{enumerate}
\end{theorem}
\textbf{Proof:}

Let us fix $a \in H$. We will show that $V^n a$ is a Cauchy sequence. For any $n,k,m,M \in \mathbb{Z}_+$ we have
\begin{multline*}
  \|V^{n+k}a - V^n a\|\\\leq
  \|V^{n+k} (P_m a) - V^{n} (P_m a) \| + \|V^{n+k} (I-P_m) a\| +  \|V^{n} (I-P_m) a\| \\
   \leq \|P_M V^{n+k} (P_m a) - P_MV^{n} (P_m a) \| + \|(I-P_M)V^{n+k} (P_m a) \| \\
    +  \|(I-P_M)V^{n} (P_m a) \|
    + \|V^{n+k} (I-P_m) a\| +  \|V^{n} (I-P_m) a\| \\
    \leq \sum_{i=1}^M \|e_i\| \sum_{j\leq m} |V^{n+k}_{ij}-V^n_{ij}| \cdot |a_j|   + \|(I-P_M)V^{n+k} (P_m a) \| \\
    +  \|(I-P_M)V^{n} (P_m a) \| + 2 K \|(I-P_m)a\|.
\end{multline*}
Let us fix $\epsilon >0$. Then by taking $m$ big enough by Lemma~\ref{lem:gcp-prop} (condition (\ref{eq:remto0})) we obtain
that $2 K \|(I-P_m)a\| \leq \epsilon$. Let us fix such value of $m$.

Now we look at  $\|(I-P_M)V^{\ell} (P_m a) \|$ for any $\ell$. Observe that  $V^{\ell} (P_m a)$ is linear combination
of first $m$ columns in $V^{\ell}$. Indeed,
\begin{equation*}
  V^{\ell}P_m a = \sum_{j \leq m} V^{\ell}_{\ast j} a_j.
\end{equation*}
From this and assumption (\ref{eq:Vnj-approribnd}) it follows that for any $\ell$ there holds
\begin{eqnarray*}
  (I-P_M)  V^{\ell}P_m a  \subset (I-P_M) \left( \sum_{j\leq m} W^j a_j \right)
\end{eqnarray*}
The set $\sum_{j\leq m} W^j a_j$ is compact, hence from Lemma~\ref{lem:compt-gss} it follows there exists $M$, such that $\|(I-P_M) \left( \sum_{j\leq m} W^j a_j \right)\| \leq \epsilon/2$. We fix such value of $M$.

From condition (\ref{eq:dfg-Vij-conv0}) it follows that $ \sum_{i=1}^M \|e_i\| \sum_{j\leq m} |V^{n+k}_{ij}-V^n_{ij}| \cdot |a_j| \leq \epsilon$ for $n$ large enough and any $k$.

Summarizing, we obtained that $\|V^{n+k}a - V^n a\|\leq 3 \epsilon$ for $n$ large enough and therefore
$V^n a$ converges to some $Va$, which defines a linear operator $V:H \to H$ satisfying $\|V\| \leq K.$

Since $\pi_i V^n e_j=e_i V^n_{ij}$, we see that $V_{ij}=\lim_{n \to \infty} V^n_{ij}$.

Finally, from Lemma~\ref{lem:formOnH} applied
to $\pi_i V$ it follows that  $\sum_{j=1}^\infty |V_{ij}| \cdot |a_j| < \infty$.

\qed

\subsection{Block decomposition}
\label{subsec:lmgssblk-decmp}
Fix $M>1$ and let $J_1,\dots,J_m$ be a partition of the set $\{1,\dots,M\}$. We set $J_{m+1}=\{M+1,M+2,\dots\}$.
For $\ell=1,\dots,m+1$ we set $H_{<\ell>}=P_{J_\ell}H$. This defines a decomposition
\begin{equation}
 H=\bigoplus_{\ell \leq m+1} H_{\langle \ell \rangle}.   \label{eq:H-decmp}
 \end{equation}
For $x \in H$ we have $x=\sum_{\ell \leq m+1} x_{\langle \ell \rangle}$, where $x_{\langle \ell \rangle}=P_{J_{\ell}}x$.

For an operator $A \in \mbox{Lin}(H,H)$ we define its blocks according to the above decomposition by setting
\begin{equation*}
  \BD{A}{k_1}{k_2} x= P_{J_{k_1}}Ax, \quad x\in H_{\langle k_2\rangle}.
\end{equation*}
With this notation we have
\begin{equation*}
   (A x)_{\langle k \rangle} = \sum_{\ell=1}^{m+1} \BD{A}{k}{\ell} x_{\langle \ell \rangle}.
\end{equation*}


\begin{theorem}
\label{thm:gcpn-abstr}
Let $(H,\|\cdot\|)$ be \gss and fix a block decomposition of $H$, as in (\ref{eq:H-decmp}). Assume that
\begin{itemize}
\item $V^n:H_n \to H_n$, $n>M$ is a family of linear maps, such that for all $i,j \in \mathbb{Z}_+$ there exist a finite limit
$V_{ij}:=\lim_{n\to \infty} V^n_{ij}$, 
\item there exists a family $W^j\subset H$, $j\in\mathbb Z_+$ of compact sets, such that for $n>M$ there holds
 $V^n_{\ast j} \in W^j$,
\item there exists $B \in \mathbb{R}^{(m+1) \times (m+1)}$, such that
      \begin{eqnarray}
       \|\BD{V^n}{k}{\ell}\| &\leq& B_{k\ell}, \quad k,\ell \leq m+1,\ n>M.  \label{eq:dfgVijn-estm}
      \end{eqnarray}
\end{itemize}
Then  there exists a bounded linear operator $V: H \to H$ such that
\begin{enumerate}
\item  for all $a \in H$ there holds $\lim_{n \to \infty} V^n(a)=Va$,
\item  for each $i$ and all $a \in H$ there holds $\sum_{j=1}^\infty |V_{ij}| \cdot |a_j| < \infty$,
\item  $\pi_i V(a) =  \left(\sum_{j=1}^\infty  V_{ij} a_j\right) e_i$ and
\item we have the following bounds
       \begin{eqnarray}
        \|\BD{V}{k}{\ell}\| &\leq& B_{k\ell}, \quad k,\ell \leq m+1.  \label{eq:dfgVij-estm}
      \end{eqnarray}
\end{enumerate}
\end{theorem}
\noindent
\textbf{Proof:}
From (\ref{eq:dfgVijn-estm}) and condition \NC4 we have
\begin{eqnarray*}
  \|V^n a\|  \leq \sum_{k,\ell=1}^{m+1} \|\BD{V^n}{k}{\ell}\| \cdot \|a_{\langle \ell \rangle}\| \leq  \left(\sum_{k,\ell=1}^{m+1} B_{k\ell}\right) \|a\|
\end{eqnarray*}
for $a \in H$. This proves that $\|V^n\|$ are uniformly bounded. From Theorem~\ref{thm:Vn-weak-lim} it follows that for each $a \in H$ the sequence $V^n a$ converges to $Va$, $V:H \to H$ is continuous and for each $i$ there holds  $\sum_{j=1}^\infty |V_{ij}| \cdot |a_j| < \infty$.

There remains to prove (\ref{eq:dfgVij-estm}). From (\ref{eq:dfgVijn-estm}) we have
$$
\|\BD{V}{k}{\ell}\| = \sup_{\|a_\ell\|=1}\|\BD{V}{k}{\ell}a_{\langle \ell \rangle}\|\leq \sup_{n>M}\sup_{\|a_\ell\|=1}\|\BD{V^n}{k}{\ell}a_\ell\|\leq B_{k\ell},
$$
which completes the proof.

\qed

%% file: lineqestm.tex
\section{Logaritmic norms and estimates for li\-near \\ equa\-tions}
\label{sec:estmlinEq}

The goal of this section is to provide estimates for solutions to non-autonomous linear ODEs, when only bounds on the time-dependent coefficients are known. Later, in Section~\ref{sec:c1-conver}, we will use these results to provide uniform bounds on the solutions to all Galerkin projections of variational equations. 

The content of this subsection is an extension of the results from \cite{KZ}. Let $\| \cdot \|$ denote a norm on $\mathbb{R}^n$  and by the same symbol we will denote operator norm induced on $\mathbb{R}^{n \times n}$. 
\begin{definition} \cite{D58,L58}
The \emph{logarithmic norm} of an operator $A$ is defined by 
\begin{equation*}
  \mu(A)=\lim_{h \to 0^+} \frac{\|\mathrm{Id} + hA \| - 1}{h}.
\end{equation*}
\end{definition}
For the properties and usage of the logarithmic norm see \cite{KZ} and the literature given there.

\begin{definition}
Let $f :U \to \mathbb{R}^n$, where  $U \subset \mathbb{R}$ is open. An absolutely continuous function $x:[a,b] \to \mathbb{R}^n$ is a \emph{weak solution} of
\begin{equation*}
x'=f(x), \quad  x\in \mathbb{R}^n
\end{equation*}
if for some $t_0\in[a,b]$ and all $t \in [a,b]$ there holds
\begin{equation*}
	x(t) = x_0 + \int_{t_0}^t f(x(s)) ds.
\end{equation*}
\end{definition}

Consider a linear equation
\begin{equation}
	v'(t)=A(t) \cdot v(t) + b(t), \label{eq:linnonhp}
\end{equation}
where $v(t) \in \mathbb{R}^k$, $A(t) \in \mathbb{R}^{k \times k}$, $b(t) \in \mathbb{R}^k$, $A$ and $b$ are bounded and measurable.

Assume that we have a decomposition of the phase space $\mathbb{R}^k$ of the form $\mathbb{R}^k=\bigoplus_{i=1}^n
\mathbb{R}^{k_i}$, $k_i\geq 1$ and accordingly for $z \in
\mathbb{R}^k$ we will write $z=(z_{\langle1\rangle},\dots,z_{\langle n\rangle})$, where $z_{\langle i\rangle} \in
\mathbb{R}^{k_i}$. Using this notion, equation (\ref{eq:linnonhp}) can be written as
\begin{equation*}
  z'_{\langle i\rangle}(t)= \sum_{j=1}^n \BD{A}{i}{j}(t)z_{\langle j\rangle}(t) + b_{\langle i\rangle}(t), \quad i=1,\dots,n,
\end{equation*}
where $z_{\langle i\rangle}(t), b_{\langle i\rangle}(t) \in \mathbb{R}^{k_i}$, $i=1,\dots,n$ and $\BD{A}{i}{j}(t) \in
L(\mathbb{R}^{k_i},\mathbb{R}^{k_j})$, $i,j=1,\ldots,n$ are linear maps. To each block $\BD{A}{i}{j}$ of the matrix $A$ we will assign a real number $J_{ij}$ and collect them in an $n\times n$ matrix $J$. The quantity $J_{ij}$ will estimate the influence of $z_{\langle j\rangle}$ on $z'_{\langle i\rangle}$. The details will be given in the sequel.

The following lemma plays a crucial role in our considerations.
\begin{lemma}\cite[Lemma 4.1]{KZ}\label{lem:boundOnVariationalLinEq}
Let $z:[0,T] \to {\mathbb R}^k=\bigoplus_{i=1}^n \mathbb{R}^{k_i}$ be a weak solution of the equation
\begin{equation*}
  z'(t)=A(t) \cdot z(t) + b(t),
\end{equation*}
where $A:[0,T]\to {\mathbb R}^{k \times k}$ and $b:[0,T] \to {\mathbb R}^k$ are bounded and measurable.

Assume that $J:[0,T] \to\mathbb{R}^{n \times n}$ and $c:[0,T]\to\mathbb R^n$ are measurable and for all $t \in [0,T]$ there holds
\begin{eqnarray*}
 J_{ij}(t) &\geq&
  \begin{cases}
      \|\BD{A}{i}{j}(t) \|  & \text{for $i \neq j$}, \\
     \mu(\BD{A}{i}{i}(t))  & \text{for $i=j$},
  \end{cases}\\
  c_i(t)&\geq&\|b_{\langle i\rangle}(t)\|\ \text{ for $i=1,\ldots,n$}.
\end{eqnarray*}
Then, for $i=1,\ldots,n$ and $t\in[0,T]$ there holds
\begin{equation*}
  \|z_{\langle i\rangle}\|(t) \leq y_i(t),
\end{equation*}
where  $y:[0,T] \to \mathbb{R}^n$ is a weak solution of 
\begin{equation*}
   y'(t)=J(t)y(t) + c(t)
\end{equation*}
with an initial condition $y(0)$ satisfying
\begin{equation*}
	\|z_{\langle i\rangle}(0)\|\leq y_i(0),\quad i=1,\ldots,n.
\end{equation*}
\end{lemma}

Next lemma gives bounds for fundamental matrix for a non-autonomous linear equation.
\begin{lemma}
\label{lem:comp-lin-problems}
Fix	a decomposition $\mathbb R^{k}=\bigoplus_{i=1}^n \mathbb{R}^{k_i}$. Let $W:[0,T] \to \mathbb{R}^{k\times k}$ be a weak solution of the equation
\begin{equation*}
  W'(t)=A(t) \cdot W(t),
\end{equation*}
where  $A:[0,T]\to {\mathbb R}^{k \times k}$ is bounded and measurable. 

Assume $J:[0,T] \to \mathbb{R}^{n \times n}$ is a measurable matrix function satisfying the following inequalities for all $t \in [0,T]$
\begin{equation*}
 J_{ij}(t) \geq
  \begin{cases}
	\|\BD{A}{i}{j}(t) \|  & \text{for $i \neq j$}, \\
	\mu(\BD{A}{i}{i}(t))  & \text{for $i=j$}.
\end{cases}\\
\end{equation*}

Then
\begin{equation*}
	\|\BD{W}{i}{j}(t)\| \leq B_{ij}(t), \quad t \in [0,T],  \quad i,j=1,\dots,n,
\end{equation*}
where $B:[0,T] \to \mathbb{R}^{n \times n}$ is a weak solution of the equation 
\begin{equation*}
   B'(t)=J(t)B(t),\quad 
\end{equation*}
with an initial condition $B(0)$ satisfying 
\begin{equation*}
  \|\BD{W}{i}{j}(0)\| \leq  B_{ij}(0), \quad i,j=1,\dots,n.
\end{equation*}
\end{lemma}
\textbf{Proof:}
For any $z^0 \in \mathbb{R}^k$ the function $z(t)$ given by
\begin{equation*}
  z(t)=W(t)z^0
\end{equation*}
is a weak solution of equation $z'=A(t)z(t)$ with the initial condition $z(0)=W(0)z^0$. Let us fix $j \in \{1,\dots,n\}$. Then, for any $z^0$ such that $\|z^0_{\langle i\rangle}\|=\delta_{ij}$ and for all $i=1,\ldots,n$ we have 
\begin{eqnarray*}
	\|z_{\langle i\rangle }(0)\| = \|(W(0)z^0)_{\langle i\rangle }\| \leq \sum_{j=1}^n \|\BD{W}{i}{j}(0)\| \cdot \| z^0_{\langle j\rangle }\| \leq B_{ij}(0).
\end{eqnarray*}
Take  $y(t)$ solving equation $y'=Jy$ with the initial condition $y(0)=B_{\ast j}(0)$. From Lemma~\ref{lem:boundOnVariationalLinEq} with $b(t)\equiv 0$ it follows that 
\begin{equation*}
\|\BD{W}{i}{j}(t)z^0_{\langle j\rangle }\| = \|(W(t)z^0)_{\langle i\rangle }\| = \|z_{\langle i\rangle}(t)\| \leq y_i(t) = B_{ij}(t).
\end{equation*}
Since $\|z^0_{\langle j\rangle}\|=1$ is arbitrary, we obtain $\|\BD{W}{i}{j}(t)\| \leq B_{ij}(t)$.
\qed

Consider the differential equation
\begin{equation}
  x'=f(x), \quad \mbox{$f \in \mathcal C^1$}. \label{eq:odelogn}
\end{equation}
Let $\varphi(t,x_0)$  denote the solution of equation
(\ref{eq:odelogn}) with the initial condition $x(0)=x_0$. By $\|x
\|$ we denote a fixed arbitrary norm in $\mathbb{R}^n$.

From Lemma~\ref{lem:boundOnVariationalLinEq} one can easily derive the following result.
\begin{lemma}
\label{lem:estmLogN} 
Let $y:[0,T] \to \mathbb{R}^n$ be a piecewise
$\mathcal C^1$ function and $\varphi(\cdot,x_0)$ be defined for $t \in
[0,T]$. Suppose that $Z$ is a convex set such that   the following
estimates hold:
\begin{eqnarray*}
  y([0,T]), \varphi([0,T], x_0) \in Z  \\
  \mu\left(\frac{\partial f}{\partial x}(\eta)\right) \leq l,\quad \mbox{ for $\eta \in Z$} \\
  \left\| \frac{dy}{dt}(t) - f(y(t)) \right\| \leq \delta.
\end{eqnarray*}
Then for $0 \leq t \leq T$ there is
\begin{displaymath}
 \| \varphi(t,x_0) - y(t)  \| \leq e^{lt}  \|y(0) - x_0 \| + \delta \kappa_l(t),
\end{displaymath}
where
\begin{equation}\label{eq:kappa}
	\kappa_l(t) = \begin{cases}
		\frac{e^{lt} -1}{l}& \mbox{if $l \neq 0$},\\
		t & \mbox{if $l =0$}.
	\end{cases}
\end{equation}
\end{lemma}

%% file: c0conver.tex
\section{$C^0$-convergence}
\label{sec:c0-conver}


The aim of this section is to state and prove theorem about uniform convergence of local flows induced by Galerkin projections to the local semi-flow for the original infinite-dimensional system. The assumptions, called later the $\mathcal C^0$-convergence conditions, will be about local properties of the vector field, only. We will assume, that the vector field satisfies certain condition on some sets (rough enclosures) satisfying the following geometric properties.

In the context of $H$ being \gss, we define a special set of Galerkin projections as follows.
Let $\emptyset\neq J_1  \subsetneq  J_2  \subsetneq J_3 \dots  $ be a family of finite sets, such that $\bigcup_{n \in \mathbb{Z}_+} J_n = \mathbb Z_+$.
For $n \in \mathbb{Z}_+$  by $H_n$ we denote a subspace spanned by
$\{e_j\}_{j \in J_n}$. Put $P_n:=P_{J_n}$ and $Q_n=\mathrm{Id}-P_n$. By $\iota_n:H_n \to H$ we denote the embedding $H_n$ into $H$ and we define $\pi_k x = P_k x - P_{k-1}x$, i.e.
$\pi_k$ is the projection onto space spanned by $\{e_j\}_{j \in J_k \setminus J_{k-1}}$.
Observe that we abuse a bit the notation here, because in Section~\ref{sec:gss} $P_n$, $H_n$, $\pi_n$ and $\iota_n$ had slightly different meaning, however the results proven there apply also in the present context, it is enough to replace $n$-th direction spanned by $e_n$ by a finite dimensional space $P_{J_n\setminus J_{n-1}}H$ which is spanned
by $\{e_j\}_{j \in J_n \setminus J_{n-1}}$. We will call family $\{J_n\}_{n \in \mathbb{Z}_+}$ a \emph{Galerkin filtration of $H$}.  The trivial Galerkin filtration of $H$ is defined by
$J_n=\{j , j \leq n\}$.

\begin{definition}
	Assume that $(H,\|\cdot\|)$ is \gss space. We say that $W\subset H$ satisfies condition $\mathbf{S}$ if
	\begin{description}
		\item[S1:] $W$ is convex and there exists $M \geq 1$, such that $P_n(W) \subset W$
		for $n \geq M$,
		\item[S2:] $W$ is compact.
	\end{description}
\end{definition}

For a vector field $F:\dom(F)\subset H\to H$ we define its $n$-th Galerkin projection $F^n:H_n\to H_n$ by
\begin{equation}\label{eq:GalerkinODE}
u'=F^n(u) :=P_n(F(i_n(u))).
\end{equation}



\begin{definition}
  We say that $F$ is admissible, if for all $n\in\mathbb Z_+$  the function $F^n$ is defined on $H_n$ and it is $\mathcal{C}^3$ as a function on $H_n$.
\end{definition}


Now we are in the position to present our $\mathcal C^0$ convergence conditions for infinite-dimensional vector fields.
\begin{definition}
Let $(H,\|\cdot\|)$ be \gss,  $W\subset H$ and $F:\mathrm{dom}(F)\subset H \to H$. We say that $F$ satisfies condition \textbf{C} on $W$ if $F$ is admissible, $W$ satisfies condition \textbf{S} and
\begin{description}
\item[C1:] $W \subset \dom(F)$ and function $F|_W:W \to H$ is continuous;
\item[C2:] there exists  $l \in \mathbb{R}$ such that for all $n\in\mathbb Z_+$ there holds
\begin{equation*}
	\sup_{x \in P_n W}\mu\left(D F^n (x) \right) \leq l.
\end{equation*}
\end{description}
\end{definition}
The main idea behind condition \textbf{C2} is to ensure that the logarithmic norms (see Section~\ref{sec:estmlinEq}) for all Galerkin projections are uniformly bounded.

Observe that conditions \textbf{S} and \textbf{C1} imply that $F \circ P_n$ converges uniformly to $F$ on $W$.

The following theorem is a generalization of  Theorem 13 in \cite{Z}. There   it is assumed that  $W$ is a trapping region and  $H=l_2$ was used, but  the main idea of the proof is the same.
\begin{theorem}
\label{thm:limitLN}
 Let $(H,\|\cdot\|)$ be \gss and assume that $F:\mathrm{dom}(F)\subset H \to H$ satisfies condition \textbf{C} on $W\subset H$. Let $Z\subset W$ and $T>0$ be such that for all $n>M_1$ there holds
  \begin{equation}
    \varphi^n(t,x)\in W\quad \text{for } x\in P_n(Z),\ t\in[0,T], \label{eq:appriori-bnds}
  \end{equation}
  where $\varphi^n$ is a local dynamical system on $H_n$ induced by the $n$-th Galerkin projection (\ref{eq:GalerkinODE}). Then there exists a continuous function $\varphi\colon[0,T]\times Z\to W$ satisfying the following properties.
    \begin{description}
    \item[1.]{\bf Uniform convergence:} The functions $\widehat{\varphi}^n:=\iota_n\circ \varphi^n\circ P_n$ converge uniformly to $\varphi$ on $[0,T]\times Z$.
    \item[2.] {\bf Existence and uniqueness within $W$:} For all $x\in Z$ the function $u(t):=\varphi(t,x)$ is a unique solution to the initial value problem $u'=F(u)$, $u(0)=x$ and satisfying $u(t)\in W$ for $t\in[0,T]$.
    \item[3.]  {\bf Lipschitz constant:} For $x,y\in Z$ and $t\in[0,T]$ there holds
      \begin{equation*}
         \|\varphi(t,x) - \varphi(t,y) \| \leq e^{lt}\|x - y\|.
      \end{equation*}
    \item[4.] {\bf Semidynamical system.} The partial map $\varphi:[0,T] \times W \to W$ defines a semidynamical system on $W$, that is
        \begin{itemize}
          \item $\varphi$ is continuous;
          \item $\varphi(0,u)=u$;
          \item $\varphi(t,\varphi(s,u)) = \varphi(t+s,u)$
        \end{itemize}
    provided $\varphi(t+s,u)$ exists.
  \end{description}
\end{theorem}
\begin{remark}
Observe that the essential difficulty in application of the above theorem is to find set $W$ satisfying (\ref{eq:appriori-bnds}) and the convergence conditions. A systematic way to construct it is based on the \emph{isolation property} (see Section~\ref{sec:isolation}). We will prove there that it is always possible to find a-priori bounds and make one time step of a rigorous integration algorithm.
\end{remark}
\noindent
\textbf{Proof of Theorem~\ref{thm:limitLN}:}

Let us fix $k \geq n$ and set $x^n=\varphi^n(\cdot,u)$, $x^k=\varphi^k(\cdot,v)$ for some $u\in P_nZ$ and $v\in P_kZ$. We start from estimation on the difference $x^n(t)-x^k(t)$ for $t \in [0,T]$. We have
\begin{eqnarray*}
  (P_n x^k(t))'&=&P_n F(x^k(t)) \\
  &=&P_nF^n(P_nx^k(t)) + \left(P_n F(x^k(t)) - P_n F(P_nx^k(t)) \right).
\end{eqnarray*}
Hence we can treat the function $y(t)=P_n x^k(t)$ as a solution to a perturbed equation $y'(t)=P_nF^n(y(t))+\delta(t)$, where $\delta$ is uniformly bounded by
\begin{eqnarray*}
  \|\delta(t)\|=\|P_n F(x^k(t)) - P_n F(P_nx^k(t))\| \leq \max_{y \in W} \|P_n F(y) - P_nF(P_ny)\| =:\delta_n.
\end{eqnarray*}
Obviously $\delta_n \to 0$ for $n \to \infty $, because $F \circ P_n$ converges uniformly to $F$ on $W$ -- this follows immediately from compactness of $W$ (condition \textbf{S2}) and \textbf{C1}.

From \textbf{C2} there is a uniform bound $l$ on all logarithmic norms of $DF^n(W)$, $n\in\mathbb Z_+$. Hence, by Lemma \ref{lem:estmLogN} we obtain
\begin{equation}
  \|x^n(t) - P_n(x^k(t))\| \leq e^{lt}\|x^n(0) - P_nx^k(0)\| + \delta_n \kappa_l(t), \quad t \in [0,T]  \label{eq:xn-Pn}
\end{equation}
where $\kappa_l$ is defined as in (\ref{eq:kappa}).

\textbf{Convergence.} For $x\in Z$, $t\in[0,T]$ and $k\geq n$ from (\ref{eq:xn-Pn}) with $u=P_nx$ and $v=P_kx$ we have
\begin{multline*}
 \|\widehat\varphi^n(t,x) - \widehat\varphi^k(t,x)\| \leq  \\
 \|\varphi^n(t,P_nx) - P_n(\varphi^k(t,P_kx))\| + \|(I-P_n) \varphi^k(t,P_kx)\| \leq \\
    \delta_n \kappa_l(t) + \|(I-P_n) x_k(t)\| \leq
     \delta_n \kappa_l(t) + \|(I-P_n)W\|.
\end{multline*}
From Lemma~\ref{lem:compt-gss} we get $\lim_{n\to\infty}\|(I-P_n)W\| = 0$. The function $\kappa_l$ is bounded on $[0,T]$ and given that $\lim_{n\to\infty}\delta_n=0$ we see that $\widehat\varphi^n$ is a Cauchy sequence in $\mathcal{C}([0,T]\times Z,H)$. Hence, it converges uniformly to a continuous function $\varphi:[0,T]\times Z\to W$.

\textbf{Existence.}
The function $\widehat\varphi^n$ satisfies the following integral equation for $x\in Z$
\begin{equation*}
	\widehat \varphi^n(t,x)=\widehat \varphi^n(0,x) + \int_0^t \iota_n P_n F(\widehat\varphi^n(s,x))ds.
\end{equation*}
From the uniform continuity of $F$ on $W$ and already established convergence of $\widehat \varphi^n(t,x)$ it follows that after passing to the limit $n \to \infty$ there holds
\begin{equation*}
	\varphi(t,x)=\varphi(0,x) + \int_0^t F(\varphi(s,x))ds,
\end{equation*}
which implies that for $x\in Z$ the function $u(t) = \varphi(t,x)$ is a solution of $u'=F(u)$ with $\varphi(0,x)=x$.

\textbf{Uniqueness.}
Let $u:[0,T] \to W$ be a solution of $u'=F(u)$ with initial condition $u(0)=x$. We will show that the sequence $\widehat \varphi^n(\cdot,x)$ converges uniformly to $u$. Applying Lemma~\ref{lem:estmLogN} to $n$-th Galerkin projection and the function $P_n u(t)$ we obtain an estimate
\begin{equation*}
  \|\widehat \varphi^n(t,x) - P_n(u(t))\| \leq  \delta_n \kappa_l(t).
\end{equation*}
Since the tail $\|(I-P_n)u(t)\| \leq \|(I-P_n)W\|$ is by Lemma~\ref{lem:compt-gss} uniformly converging to zero as $n\to \infty $, we see that $\widehat\varphi^n(\cdot,x) \to u$ uniformly on $[0,T]$.

\textbf{Lipschitz constant on $W$}. From Lemma~\ref{lem:estmLogN} applied to $n$-dimensional Galerkin projection for different initial conditions we obtain
\begin{equation*}
\|\widehat\varphi^n(t,x)-\widehat \varphi^n(t,y)\| =  \|\varphi^n(t,P_nx)-\varphi^n(t,P_ny)\| \leq e^{lt}\|P_n x - P_n y\|,
\end{equation*}
which after passing to the limit $n\to\infty$ gives
\begin{equation*}
  \|\varphi(t,x)-\varphi(t,y)\| \leq e^{lt}\|x - y\|.
\end{equation*}

Assertion 4 follows easily from the previous ones. \qed

%% file: c1conver.tex
\section{$\mathcal C^1$-convergence of Galerkin projections}
\label{sec:c1-conver}

In Section~\ref{sec:c0-conver} we gave (Theorem~\ref{thm:limitLN}) sufficient and verifiable by means of rigorous numerics conditions that guarantee uniform convergence of a sequence of flows induced by Galerkin projections to a solution of the underlying infinite dimensional system. In this section we would like to address similar question for associated variational equations. Before we state the assumptions about the vector field we need to introduce some notion.

\subsection{Some remarks about $C^1$ functions on compact subsets of \gss spaces}

Let $f:H \supset \dom (f) \to \mathbb{R}$. In the applications we keep in mind, the function $f$ will be a component of a vector field, which might be not  Frech\'et differentiable on $H$. On the other hand, assuming that the restrictions $f^n:=f\circ i_n: H_n \to \mathbb{R}$ are $\mathcal C^1$-smooth for all $n$, the derivative $Df^n(z)$ can be seen as a linear form on $H_n$ represented by a row-vector   $\left(\frac{\partial f}{\partial x_1}(z),\dots,  \frac{\partial f}{\partial x_n}(z)\right)$. Formally $Df^n(z) \in H_n^*$ for $z \in H_n$. We can embed it to $H^*$ by setting $Df^n(z) \circ P_n $ for $z \in H_n$.

Here a natural question arises: does the infinite row-vector  \newline $\left(\frac{\partial f}{\partial x_1}(z), \frac{\partial f}{\partial x_2}(z),\dots\right)$ define a $1$-form on $H$? If so, can we treat it as the derivative of $f$? For this purpose we introduce the notion
\begin{equation}\label{eq:rowderivative}
	\widetilde{D f}(z) := \lim_{n \to \infty} Df^n(P_n z) \circ P_n.
\end{equation}
and prove the following lemma.

\begin{lemma}
	\label{lem:df-form}
	Let $H$ be \gss, $W\subset H$ and let $f:H \supset \dom (f) \to \mathbb{R}$.
	Assume that
	\begin{itemize}
		\item $W$ satisfies condition \textbf{S};
		\item $f|_W$ is continuous;
		\item $f^n$ is $\mathcal C^1$ on $H_n$;
		\item for all $z \in W$ the limit (\ref{eq:rowderivative}) exists and the function $ \widetilde{D f}: W \to H^*$ is continuous.
	\end{itemize}
	Then for every  $z,w \in W$ there holds
	\begin{eqnarray}
		f(z) - f(w) &=& \left(\int_{0}^1 \widetilde{Df}(t z + (1-t)w)dt\right) (z-w) \nonumber \\
		&=&\sum_{j\in\mathbb Z_+}  \int_0^1 \widetilde{D_j f}(t z + (1-t)w)dt \cdot (z_j - w_j),    \label{eq:diff-f-onW}
	\end{eqnarray}
	where $\widetilde{D_j f}(z) := \widetilde{D f}(z)(e_j)$. Moreover, if for some $z \in W$ and $j$ there exists $\delta>0$, such that $ z+[-\delta,\delta]e_j \in \dom (f)$, then   $\frac{\partial f}{\partial x_j}(z)$  exists and
	\begin{equation}
		\widetilde{D_j f}(z)= \frac{\partial f}{\partial x_j}(z).  \label{eq:partderF}
	\end{equation}
\end{lemma}
\textbf{Proof:}
Observe that for any $v,z \in H$ and $n,k \in \mathbb{Z}_+$ there holds
\begin{eqnarray*}
	Df^{n+k}(P_{n+k}(P_nz)) (P_{n+k}(P_nv))&=&Df^{n+k}(P_nz) (P_nv)\\
	&=& Df^{n}(P_nz) (P_nv).
\end{eqnarray*}
Hence
\begin{equation}
	\widetilde{D f}(P_nz)(P_n v)= Df^n(P_nz) (P_nv).  \label{eq:wtDf=Df-onP}
\end{equation}

For any $n \geq M$ and  $z,w \in W$ using (\ref{eq:wtDf=Df-onP}) we have
\begin{eqnarray*}
	f(P_n z) - f(P_n w)=\sum_{j \leq n} \int_0^1 \frac{\partial f}{\partial x_j}(P_n (t z + (1-t)w))dt \cdot (z_j - w_j) \\
	= \int_{0}^1 Df^n(P_n (t z + (1-t)w))dt \cdot P_n(z-w)  \\
	= \int_{0}^1 \widetilde{Df}(P_n (t z + (1-t)w))dt \cdot P_n(z-w).
\end{eqnarray*}
Now we pass to the limit $n \to \infty$ in the above equation. From continuity of $f$ it follows that  $f(P_n z) - f(P_n w) \to f(z) - f(w)$ and from \NC2 we obtain $ P_n(z-w)  \to (z-w)$. From continuity of $\widetilde{Df}$ on $W$, which is compact,  follows its uniform continuity. Moreover, from Lemma~\ref{lem:compt-gss} we get that $P_n$ converges uniformly to the identity on $W$. Therefore the integral converges to  $\int_{0}^1 \widetilde{Df}(t z + (1-t)w)dt$ and we obtain
\begin{eqnarray*}
	f(z) - f(w)= \left(\int_{0}^1 \widetilde{Df}(t z + (1-t)w)dt\right) (z-w).
\end{eqnarray*}
Observe that $\int_{0}^1 \widetilde{Df}(t z + (1-t)w)dt \in H^*$, because it is bounded \\
by $\max_{z \in W} \|\widetilde{Df}(z)\|$.  Therefore from Lemma~\ref{lem:formOnH} we obtain assertion (\ref{eq:diff-f-onW}).

Assume that $z \in W$ and $z+he_j \in W$. Then from (\ref{eq:diff-f-onW}) we obtain
\begin{eqnarray*}
	f(z+he_j) -f(z)= h \left(\int_{0}^1 \widetilde{Df}(z + th e_j)dt\right)(e_j) \\
	= h \int_{0}^1 \widetilde{D_jf}(z + th e_j)dt = h \widetilde{D_jf}(z) + h\int_0^1\left( \widetilde{D_jf}(z + th e_j) - \widetilde{D_jf}(z) \right)dt.
\end{eqnarray*}
From the continuity of $ \widetilde{D_jf}$ it follows that the integral converges to $0$ if $h \to 0$. Therefore we obtain (\ref{eq:partderF}).
\qed

\subsection{Variational equations and Galerkin projections}
%
%
%
%

Let $\varphi^n$ be a local flow induced by (\ref{eq:GalerkinODE}) and consider the variational matrix for $\varphi^n$ given by
\begin{equation*}
  V^n_{ij}(t,u)= \frac{\partial \varphi^n_i}{\partial u_j} (t,u).
\end{equation*}
Since $F$ is admissible, each column in $V^{n}_{ij}$ evolves separately and satisfies the system of variational equations
\begin{eqnarray}
  \frac{d}{dt} x^n &=& F^n(x^n), \quad x^n \in H_n,  \label{eq:sysVar1} \\
  \frac{d}{dt} C^n &=& D F^n(x^n) C^n, \quad C^n \in H_n. \label{eq:sysVar2}
\end{eqnarray}
The above system has $H_n \times H_n \subset H \times H$ as the phase-space.

Following discussion at the beginning of this section, instead of $DF_i$ which may not exists, we will use the notion $\widetilde{DF_i}$ and define a vector field on $H \times H$ as
\begin{equation}
  \frac{d}{dt} (x,C)=F_V(x,C):=(F(x),\{\widetilde{DF_i}(x)C\}_{i\in\mathbb Z_+}).  \label{eq:sysVarH}
\end{equation}

To bring the variational system to the setting discussed in Section~\ref{sec:c0-conver} we proceeded as follows. On $H \times H$ we use the norm $\|(x,v)\|=\max(\|x\|,\|v\|)$. This makes $H\times H$ a \gss space and its basis is parameterized by $\mathbb{Z}_+ \times \mathbb{Z}_+$. Given $\{J_n\}$ Galerkin filtration of $H$, we define $\{J'_n\}$ Galerkin filtration of $H \times H$ by $J'_n=J_n \times J_n$.  Then accoridingly  Galerkin projection $P_n$ on $H\times H$ is defined by  $P_n(x,C)=(P_nx,P_nC)$.

For the remainder of this section we adopt notation used in Section~\ref{sec:c0-conver} and Theorem~\ref{thm:limitLN}. In particular we have a Galerkin filtration of $H$ and
the induced Galerkin filtration of $H \times H$, so that $P_n$ will always refer to these filtrations.

\begin{definition}
\label{def:condV}
Let $(H,\|\cdot\|)$ be \gss,  $W\times W_V\subset H\times H$ and $F:\mathrm{dom}(F)\subset H \to H$. We say that $F$ satisfies condition \textbf{VC} on $W\times W_V$ if $F$ is admissible, $W\times W_V$ satisfies condition \textbf{S} and
\begin{itemize}
	\item[\bf{VC1:}] the function $F_V$ is continuous on $W\times W_V$;
	\item[\bf{VC2:}] there exists a constant A, such that
	$$
		\sup_{(x,v)\in P_n(W\times W_V)}\mu(DP_nF_V|_{H_n\times H_n})\leq A.
	$$
\end{itemize}
\end{definition}

The conditions \textbf{VC1} and \textbf{VC2} are nothing more than conditions \textbf{C1} and \textbf{C2}  for the vector field $F_V$. Applying arguments from Section~\ref{sec:c0-conver} to $F_V$ we will obtain the existence of solution to variational equation for each column (directional derivative), independently. Obtaining $\mathcal C^1$-like information about full system needs some additional reasoning.

The following theorem is a generalization of Theorem 2 in \cite{Z1}.
\begin{theorem}
\label{thm:c1conver} The same assumptions as in Theorem~\ref{thm:limitLN}. Let $\{W_{V_j}\}_{j \in \mathbb{Z}_+}\subset H$ be a family of sets such that for each $j \in \mathbb{Z}_+$ $F$ satisfies \VL on $W \times W_{V_j}$ and for any $n >j$ the solution of the variational problem (\ref{eq:sysVar1}--\ref{eq:sysVar2}) with initial conditions $x^n(0) \in Z$ and $C^n(0)=e_j$ satisfies
\begin{equation}
  C^n(t) \in W_{V_j}, \quad t \in [0,T]. \label{eq:var-a-priori-bnds}
\end{equation}
Let $l$ be a constant from condition {\rm \textbf{C2}} (uniform bound on logarithmic norm of $DF$ on $W$).

Then there exists $V:[0,T] \times Z  \to \mathrm{Lin}(H,H)$, such that $V_{ij}:[0,T]\times Z \to {\mathbb R}$ for $i,j \in \mathbb{Z}$ are continuous and the following properties are satisfied.
\begin{description}
\item[Convergence:]
For each $j$ the function $\widehat{V}^n_{\ast j}(t,x):= \iota_n V^n_{\ast j}(t,P_n x)$ converges to $V_{\ast j}(t,x)$ uniformly on $[0,T] \times Z$, and
\begin{eqnarray*}
  \|V(t,x)\| &\leq& e^{l t}, \quad (t,x) \in [0,T] \times Z,\\
   V(t,x)e_j &=& \sum_{i\in\mathbb Z_+} V_{ij}(t,x)e_i,
\end{eqnarray*}
and for every $a \in H$ the map $[0,T] \times Z \ni (t,x) \mapsto V(t,x)a$ is continuous.

\item[Smoothness:] For any $x,y \in Z$ and any $t\in[0,T]$ there holds
\begin{equation}
  \varphi(t,x) - \varphi(t,y) =  \int^1_0 V(t,y+s(x-y))ds \cdot (x - y),
    \label{eq:intpar}
\end{equation}
and for every $j$ the partial derivative of the flow $\frac{\partial \varphi}{\partial x_j}(t,x)$ exists and
\begin{equation}
  \frac{\partial \varphi}{\partial x_j}(t,x)=V_{\ast j}(t,x).  \label{eq:dfiduj=Vij}
\end{equation}
\item[Equation for $V$:]  $V(t,u)$ satisfies the following variational equation
  \begin{equation}
     \frac{d V_{*j}}{dt}(t,u) = \sum_i e_i \sum_k \frac{\partial F_i}{\partial u_k}(\varphi(t,u))
           V_{kj}(t,u), \label{eq:varinfdim}
  \end{equation}
 with the initial condition $V(0)=\mathrm{Id}$ in the following sense:  for each $j$ the derivative  $\frac{d V_{*j}}{dt}(t,u)$ exists, the series on r.h.s. of (\ref{eq:varinfdim}) converges uniformly
 on $[0,T] \times Z$ and equation (\ref{eq:varinfdim}) is
 satisfied.
\end{description}
\end{theorem}
\noindent
\textbf{Proof:}
We apply Theorem~\ref{thm:limitLN} to the variational system (\ref{eq:sysVarH}) for $j$-th column (i.e. with initial condition $C(0)=e_j$) separately. From $\VL$ it follows that all its assumptions are satisfied for system (\ref{eq:sysVarH}) with $Z_V=Z \times \{e_j\}$, $W_V=W \times W_{V_j}$ and with $l$ from condition $\VL\mathbf{2}$. Therefore, there exists a family of continuous functions $V_{\ast j}:[0,\infty)\times Z \to H$ for $j \in \mathbb{Z}_+$, such that for each $j$ the functions $\widehat V^{n}_{\ast j}$ converge to $V_{\ast j}$ uniformly on $[0,T] \times Z$ and $(x(t),V_{\ast j}(t,x))$ satisfies variational equation (\ref{eq:sysVarH}).

Now, we will prove (\ref{eq:intpar}) and (\ref{eq:dfiduj=Vij}). Let us fix $t \in [0,T]$ and $x,y\in Z$. For any $n$ we have
\begin{equation*}
  \varphi^n(t,P_n x) - \varphi^n(t,P_n y) =   \int^1_0 V^n(t,P_n y+s(P_n x- P_n y))ds \cdot P_n(x-y)
\end{equation*}

We will pass to the limit $n \to \infty$ in above equation. The limit of l.h.s.
is settled by Theorem \ref{thm:limitLN}.

Now we consider the r.h.s.  We will use Theorem~\ref{thm:Vn-weak-lim} to pass to the limit $n \to \infty$. We already have proved the convergence of $V^n_{ij}$. There remains to show, that $\|V^n\|$ are uniformly bounded.

From condition \textbf{VC2} and Lemmas~\ref{lem:comp-lin-problems} and~\ref{lem:estmLogN} it follows that for all $n\in\mathbb Z_+$,
\begin{equation}
  \|V^n(t,z)\| \leq e^{lt}, \qquad \forall n, (t,z) \in [0,T] \times P_nZ. \label{eq:VnlogN}
\end{equation}

Let $\overline{V}^n(t,[x,y])=\int^1_0 V^n(t,y+s(x- y))ds$.  From (\ref{eq:VnlogN})  we have
\begin{equation*}
     \|\overline{V}^n(t,[P_nx,P_ny])\| \leq e^{lt}, \quad t\in [0,T].
\end{equation*}
From convexity of $W_{V_j}$ it follows that for all $j$ and for $x,y \in Z$ there holds
\begin{equation*}
  \widehat V_{\ast j}^n(t,x),\ \iota_n \overline{V}^n(t,[P_nx,P_ny])_{\ast j} \in W_{V_j}.
\end{equation*}
Since $\widehat V^n_{ij}(t,\cdot)$ are continuous on $Z$ for all $i,j$ and converge uniformly to $V_{ij}(t,\cdot)$ we obtain
\begin{equation*}
  \int^1_0 \widehat V^n_{ij}(t,y+s(x- y))ds \to \int^1_0 V_{ij}(t,y+s(x- y))ds.
\end{equation*}
Now we use Theorem~\ref{thm:Vn-weak-lim} to conclude that $V(t,x) \in \mbox{Lin}(H,H)$, $\int^1_0 V(t,y+s(x- y))ds \in \mbox{Lin}(H,H)$ and
\begin{equation*}
  \lim_{n \to \infty} \int^1_0 V^n(t,y+s(x- y))ds \cdot (x-y) = \int^1_0 V(t,y+s(x- y))ds \cdot (x-y).
\end{equation*}
Gathering this with
\begin{eqnarray*}
\lim_{n\to\infty}  \varphi^n(t,P_n x) -  \varphi^n(t,P_n y)&=&\varphi(t,x) - \varphi(t,y),\\
\lim_{n\to\infty} P_n (x-y)&=& x-y
\end{eqnarray*}
we obtain (\ref{eq:intpar}).

Now we will show that from (\ref{eq:intpar}) we can conclude (\ref{eq:dfiduj=Vij}). Indeed, we have
\begin{eqnarray*}
  \varphi(t,x+he_j) -  \varphi(t,x) = h \int_0^1 V(t,x+she_j)ds \cdot e_j \\
  = h  V(t,x)e_j + h \int_0^1 \left( V(t,x+she_j)e_j - V(t,x)e_j \right)ds\\
  = h  V_{\ast j}(t,x) + h \int_0^1 \left( V_{\ast j}(t,x+she_j) - V_{\ast j}(t,x) \right)ds
\end{eqnarray*}
and the result follows from continuity of $V_{\ast j}(t,\cdot)$.

It remains to prove that for every $a \in H$ the map $[0,T] \times Z \ni (t,z) \mapsto V(t,z)a \in H$ is continuous. We know that for every $j$ the map $(t,x) \mapsto V_{\ast j}(t,z)$ is continuous. Let $(t_1,z_1), (t_2,z_2) \in [0,T] \times Z$. Then  for any $n$ holds
\begin{eqnarray*}
  \| V(t_1,z_1)a - V(t_2,z_2)a \| \leq \| V(t_1,z_1)P_n a - V(t_2,z_2)P_n a \| \\
   +  \| V(t_1,z_1)(I-P_n) a \| + \| V(t_2,z_2)(I-P_n) a \| \leq  \\
       \sum_{j \in J_n} \|V(t_1,z_1)_{\ast j} - V(t_2,z_2)_{\ast j}\| \cdot |a_j| + (e^{t_1 l} + e^{t_2 l}) \|(I-P_n)a\|
\end{eqnarray*}
For a given $\epsilon>0$ we take $n$ so large that $ \|(I-P_n)a\|e^{\max(0,l)T} < \epsilon/2$. We fix such $n$. Then from the continuity of $V_{\ast j}$ for $j\in J_n$
we can make the sum less that $\epsilon/2$ when $(t_2,z_2) \to (t_1,z_1)$. This finishes the proof.
\qed

\subsection{Block decomposition}
\label{subsec:blockdecmp}
Now, we would like to adapt Lemma~\ref{lem:comp-lin-problems} to the variational equation for dissipative PDEs and its solution obtained as the limit of solutions of variational equation of Galerkin projections. We will adopt the notation related to a decomposition of $H$ from Section~\ref{subsec:lmgssblk-decmp}.

\begin{theorem}
\label{thm:lineqestm-inf}
 Same assumptions and notation as in Theorem~\ref{thm:c1conver}.  Consider a decomposition of $H$ given by (\ref{eq:H-decmp}).

 Assume that matrix $J \in \mathbb{R}^{(m+1) \times (m+1)}$ satisfies
\begin{equation*}
 J_{k\ell} \geq
  \begin{cases}
     \sup_{n>M} \sup_{w \in W} \left\|\frac{\partial F^n_{\langle k \rangle}}{\partial u_{\langle \ell \rangle}}(w) \right\|,  & \text{for $k \neq \ell$, $l,\ell \leq m$+1}, \\
     \sup_{n>M} \sup_{w \in W}  \mu\left(\frac{ \partial F^n_{\langle k \rangle} }{\partial u_{\langle k \rangle}}(w) \right),     & \text{for $k=\ell\leq m+1$}
  \end{cases}
\end{equation*}

Then for any  $x \in Z$ and  $t \in [0,T]$
holds
\begin{eqnarray}
  \|V_{\langle k \rangle \langle \ell \rangle}(t,x)\| &\leq&  \left(e^{Jt}\right)_{k\ell},  \quad k,\ell=1,\dots,m+1.   \label{eq:Vij-estm}
\end{eqnarray}
\end{theorem}
\textbf{Proof:}
Let us take any $n > M$. We consider the decomposition of $H_n$
\begin{equation*}
H_n=\bigoplus_{k \leq m} H_{\langle k \rangle}   \oplus Y_n, \quad Y_n=P_n  H_{\langle m+1 \rangle}.
\end{equation*}

From Lemma~\ref{lem:comp-lin-problems} applied to $n$-th Galerkin projection with matrix $J$ we obtain that conditions (\ref{eq:Vij-estm}) are satisfied by  $V^n(t,P_n x)$ for $x \in Z$ and $t \in [0,T]$

Now we want to pass to the limit  $n \to \infty$.   We know from Theorem~\ref{thm:c1conver} that
\begin{equation*}
V^n_{ij}(t,P_nx) \to V_{ij}(t,x), \quad \forall (i,j)\in \mathbb{Z}_+^2  \quad \forall t \in [0,T], x\in Z. 
\end{equation*}
The result now follows from Theorem~\ref{thm:gcpn-abstr}.
\qed

%% file: isolation.tex
\section{Isolation property}
\label{sec:isolation}

In convergence theorems discussed in Section~\ref{sec:c0-conver} and Section~\ref{sec:c1-conver} the crucial assumption was the existence of uniform a-priori bounds $W$ and $W_V$. It turns out that such sets can be found in algorithmic way using the isolation property of the vector field discussed in Section~\ref{sec:dissipativePDEs}. The existence of such a-priori bounds will be later used in construction of an algorithm for rigorous integration of the flow and associated variational equations.


We adopt a notation used in Sections~\ref{sec:c0-conver} and Section~\ref{sec:c1-conver}. In particular we have some Galerkin filtration, which is used to define $\pi_k$. Thus, for $x \in H$,
$\pi_k x$ is the $k$-th component of $x$ and it is a vector of dimension $\# (J_k \setminus J_{k-1})$.

In what follows for  $W \subset H$ we set $W_k=\pi_k (W)$. We define the isolation property in the following way.
\begin{definition}
	Let $W\subset H$ be a set satisfying condition {\bf S} in a \gss space $H$. We say that  vector field $F:H\to H$ satisfies the \emph{isolation property} on the set $W$ if $F$ satisfies condition \textbf{C} on $W$ and there exists $K_0\in\mathbb N$ such that
	\begin{eqnarray*}
		\forall_{k\geq K_0}\, \exists_{T_k>0}\, \forall_{t\in(0,T_k]}\forall_{u\in W}
		\left(  u_k\in\mathrm{bd} W_k\Longrightarrow u_k + tF_k(u) \in\mathrm{int}W_k\right).
	\end{eqnarray*}	
\end{definition}
Geometrically the above condition means that the vector field $F$ is pointing inwards each component $W_k$ on all far coordinates of $W$. In particular, each component $W_k$, $k\geq K_0$ has nonempty interior.
\begin{lemma}\label{lem:isolation-for-projections}
	If $F:H\to H$ satisfies  the isolation property on $W\subset H$ then there is $N_0$ such that for any Galerkin projection $n>N_0$ there holds
	\begin{eqnarray*}
	\forall_{N_0\leq k\leq n}\, \exists_{T_k>0}\, \forall_{t\in(0,T_k]}\forall_{u\in P_nW}
	\left(  u_k\in\mathrm{bd} W_k\Longrightarrow u_k + tF_k^n(u) \in\mathrm{int}W_k\right).
\end{eqnarray*}	
\end{lemma}
\textbf{Proof:}
From \textbf{S1} there is $M\geq 1$ such that $P_nW\subset W$ for $n\geq M$. Let $K_0$ be the constant from  the isolation property for $F$. Put $N_0=\max\{K_0,M\}$. Now, let us fix $n>N_0$, $N_0\leq k\leq n$ and $u\in P_nW$ such that $u_k\in\mathrm{bd}W_k$. Let $T_k$ be the same constant as in the isolation property for $F$ and index $k$.

Since $k\leq n$ and $u\in P_nW\subset W$ we have $F^n_k(u) =F_k(u)$ and thus
\begin{equation*}
	u_k+tF_k^n(u) = u_k+tF_k(u)\in \mathrm{int}W_k
\end{equation*}
for $t\in(0,T_k]$.
\qed

\begin{definition}
	We say that a vector field $F:H\to H$ is isolating on the family of sets $\mathcal W\subset 2^H$ if $F$ satisfies the isolation property on each set $W\in\mathcal W$.
\end{definition}
Notice, that the constant $K_0$ for each set in the family $\mathcal W$ can be different.

The next theorem addresses the question about the existence of a-priori bounds, as in Theorem~\ref{thm:limitLN}.

\begin{theorem}\label{thm:isolationGivesAPB}
	Assume the vector field $F:H\to H$ is isolating on the family of sets $\mathcal W$ satisfying the following properties:
	\begin{itemize}
		\item there exists $E\in \mathcal W$, such that $0\in\mathrm{int}E_k$ for all $k\in\mathbb Z_+$;
		\item $\mathcal W$ is closed with respect to the addition
		\begin{eqnarray*}
			W_1,W_2\in \mathcal W &\Longrightarrow&  W_1+W_2:=\{w_1+w_2: w_1\in W_1,\, w_2\in W_2 \}\in\mathcal W.
		\end{eqnarray*}
	\end{itemize}		
	Then for any $Z\in\mathcal W$ there exists $T>0$ such that the assumptions of Theorem~\ref{thm:limitLN} are satisfied for the set $W=Z+E$.
\end{theorem}
\textbf{Proof:}
We have to prove (\ref{eq:appriori-bnds}) with $W=Z+E$. From Lemma~\ref{lem:isolation-for-projections} there is $N_0$ such that for $n>N_0$ and $N_0\leq k\leq n$ the vector field $F^n$ is pointing inwards $W_k$ on the boundary of $W_k$. Thus, for $z\in P_nZ$ the only way a trajectory $z^n(t)$ of the $n$-th Galerkin projection may escape the set $P_nW$ is through one of the leading coordinate $z_k$, $k< N_0$.

Since $W$ is compact, $F$ is continuous on $W$ and $F^n$ converge uniformly to $F$ on $W$ we have
\begin{equation*}
	\sup_{n>N_0}\sup_{k< N_0}\sup_{u\in W}\|F_k^n(u)\| < \infty.
\end{equation*}
Hence, there is $T>0$ such that for all $n>N_0$ and all $k< N_0$ (finite number of coefficients) there holds
\begin{equation*}
	[0,T]F^n_k(P_nW)\subset \mathrm{int}E_k.
\end{equation*}
From this we get
\begin{equation*}
	Z_k+[0,T]F^n_k(P_nW)\subset \mathrm{int}W_k, \quad k<N_0,
\end{equation*}
which proves that for $z\in P_nZ$ the solution of the $n$-th Galerkin projection $z^n(t)$ exists for $t\in[0,T]$ and $z^n([0,T])\subset P_nW$.
\qed

In the applications we keep in mind, the family $\mathcal W$ may consists of sets with polynomial, exponential or mixed decay of far coefficients, which are clearly closed with respect to addition and the existence of $E$ with $0\in E_k$ is also satisfied.

In a similar way we can address the question regarding the existence of a-priori bounds for variational system needed  in Theorem~\ref{thm:c1conver}.
\begin{theorem}
	Assume the vector field $F:H\to H$ is isolating on the family of sets $\mathcal W$ satisfying the following properties:
	\begin{itemize}
		\item there exists $E\in \mathcal W$, such that $0\in\mathrm{int}E_k$ for all $k\in\mathbb Z_+$;
		\item $\mathcal W$ is closed with respect to the addition
		\begin{eqnarray*}
			W_1,W_2\in \mathcal W &\Longrightarrow&  W_1+W_2:=\{w_1+w_2: w_1\in W_1,\, w_2\in W_2 \}\in\mathcal W;
		\end{eqnarray*}
		\item for all $W\in\mathcal W$ the vector field $F_V$ is isolating on $W\times E$ and satisfies condition \VL
	\end{itemize}		
	Then for any $Z\in\mathcal W$ there exists $T>0$ such that the assumptions of Theorem~\ref{thm:c1conver} are satisfied with $W=Z+E$ and $W_{V_j}=C_jE$, for some $C_j>0$, $j\in\mathbb Z_+$.
\end{theorem}

\textbf{Proof:}
From Theorem~\ref{thm:isolationGivesAPB} there is $T>0$ and $N_0$ such that the assumptions of Theorem~\ref{thm:limitLN} are satisfied, that is (\ref{eq:appriori-bnds}) holds true with $W=Z+E$.

From Lemma~\ref{lem:isolation-for-projections} and since $F_V$ is isolating on $W\times E$ it follows that there is, possibly larger $N_1$ such that the vector field $F_V^n$ is pointing inwards $W\times E$ for all $n>N_1$ and $N_1<k\leq n$. Take $Z_V:=\frac{1}{2}E$. Reasoning as in the proof of Theorem~\ref{thm:isolationGivesAPB} and shrinking $T$ if necessary we obtain that for all $n>N_1$ and $t\in[0,T]$ there holds
\begin{eqnarray*}
	x^n(t)\in W&\text{for}& x^n(0)\in P_nZ,\\
	C^n(t)\in E&\text{for}& C^n(0)\in P_nZ_V.
\end{eqnarray*}
Let us fix $j\in \mathbb Z_+$. Since each component of $E$ is a convex set containing zero we can find a constant $C_j$ such that $e_j\in C_jZ_V$. Put $W_{V_j}=C_jE$. Now, due to linearity of the variational equation, for $n>N_1$ and $j\leq n$ we have
\begin{eqnarray*}
	C^n(t)\in W_{V_j}&\text{for}& C^n(0)\in P_n(C_jZ_V).
\end{eqnarray*}
In particular, if $C^n(0)=e_j$ then $C^n(t)\in W_{V_j}$ for $t\in[0,T]$ and thus (\ref{eq:var-a-priori-bnds}) is satisfied.

\qed 

%% file: convercond.tex
\section{Dissipative PDEs on the torus }
\label{sec:dissipativePDEs}

The aim of this section is to prove that the framework introduced in Section~\ref{sec:c0-conver} and Section~\ref{sec:c1-conver} is applicable to a certain class of PDEs. We start with technical estimates and then we show that conditions \textbf{S}, \textbf{C} and \textbf{VC} as well as the isolation property are satisfied for this class on certain families of sets.

Consider
\begin{equation}
  u_t = L u + N\left(u,Du,\dots,D^ru\right),  \label{eq:genpde}
\end{equation}
where $u \in \mathbb{R}^n$,  $x \in \mathbb{T}^d=\left(\mathbb{R}\mod 2\pi\right)^d$, $L$ is a linear operator, $N$ is a
polynomial and by $D^s u$ we denote $s^{\text{th}}$ order derivative of
$u$, i.e. the collection of all spatial (i.e. with respect to variable $x$) partial derivatives of $u$ of
order $s$. 

We require, that the operator $L$ is diagonal in the Fourier basis
$\{e^{ikx}\}_{k \in \mathbb{Z}^d}$,
\begin{equation*}
  L e^{ikx}= -\lambda_k e^{ikx},
\end{equation*}
with
\begin{eqnarray}
 L_* |k|^p &\leq& \lambda_k \leq  L^* |k|^p, \qquad \text{for all $|k| > K$ and  $K,L_*,L^* \geq 0$}, \label{eq:lambdak} \\
    p &>& r.  \label{eq:p>r}
\end{eqnarray}

If $a(t,x)$ is a sufficiently regular solution of
(\ref{eq:genpde}), then  we can expand it in Fourier series
$a(t,x)=\sum_{k \in \mathbb{Z}^d} a_k(t)e^{\mathrm{i}k\cdot x}$, $a_k \in \mathbb{C}^n$ to obtain
an infinite ladder of ordinary differential equations for the
coefficients $a_k$
\begin{equation}
  \frac{d a_k}{dt}=-\lambda_k a_k + N_k(a), \quad k \in \mathbb{Z}^d,  \label{eq:fugenpde}
\end{equation}
where $N_k(a)$ is $k$-th Fourier coefficient of function
$N(a,Da,\dots,D^ra)$.

Observe that $a_k$'s might not be independent variables. For example, assumption $a(t,x) \in \mathbb{R}^n$ forces the following \emph{reality condition}
\begin{equation}
  a_{-k}=\overline{a}_{k}.  \label{eq:reality}
\end{equation}
In such situation we have to consider the subspace defined by condition (\ref{eq:reality}). This subspace is invariant for all Galerkin projections of (\ref{eq:genpde}) onto subspaces containing both $a_k$ and $a_{-k}$. 

Other constraints like oddness or evenness of $a(t,x)$ may cause the
change of set of basic functions to $\sin(kx)$ or $\cos(kx)$ or combinations thereof.

In any case we will have $u(t,x)=\sum_{k \in I} a_k(t) e_k(x)$, where $I$ is a countable set and $a_k \in \mathbb{C}^n$ or $a_k \in \mathbb{R}^n$, hence our sequence space $H$ will be build
of real and imaginary parts of components of $a_k$ which  after choosing a suitable norm fits in the framework discussed in previous sections.

\subsection{Some examples}

\subsubsection{Kuramoto-Sivashinsky equation}
Let us consider the one-dimensional Kuramoto-Sivashinsky PDE \cite{KT,S}, which is given by
\begin{equation}\label{eq:KS}
u_t = -\nu u_{xxxx} - u_{xx} + (u^2)_x, \qquad \nu>0,
\end{equation}
where $x \in \mathbb{R}$, $u(t,x) \in \mathbb{R}$ and we impose periodic boundary conditions

In the Fourier basis we obtain  the following ladder of ordinary differential equations for complex coefficients $a_k$
\begin{equation}
  \dot{a}_k = k^2(1-\nu k^2) a_k + \mathrm{i}k \sum_{\ell \in \mathbb{Z}} a_k a_{k - \ell},  \label{eq:Ksper}
\end{equation}
and the reality constraint (\ref{eq:reality}). We can either treat (\ref{eq:Ksper}) as equation acting on sequence space indexed 
by $k \in \mathbb{Z}$ and consider invariant subspace defined by reality condition (\ref{eq:reality}) with Galerkin filtration $J_n=\{ k \in \mathbb{Z}, |k| \leq n \}$
or we can eliminate variable $a_k$ for $k<0$ and rewrite the convolution term in (\ref{eq:Ksper}) in terms of $a_k$ with nonnegative $k$'s.

If for (\ref{eq:KS}) we impose odd and periodic boundary conditions (\ref{eq:KSbc})
 \begin{equation}
u(t,x)=-u(t,-x), \qquad u(t,x) = u(t,x+ 2\pi), \label{eq:KSbc}
\end{equation}
then we can represent  $u(t,x)=\sum_{k \in \mathbb{Z}_+} -2a_k(t) \sin(kx)$,
where $a_k \in \mathbb{R}$ and equation (\ref{eq:KS})
becomes \cite{CCP,ZM}
\begin{equation*}
  \frac{d a_k}{dt}=k^2(1-\nu k^2) a_k - k \sum_{n=1}^{k-1} a_n
  a_{k-n} + 2k \sum_{n=1}^{\infty} a_n  a_{n+k}, \quad k=1,2,3\dots
\end{equation*}

Observe that conditions (\ref{eq:lambdak}) and
(\ref{eq:p>r}) are satisfied for the KS equation. Namely we have
$\lambda_k \sim \nu k^4$, $p=4$, $r=1$.

\subsubsection{Navier Stokes equations on the torus}
 The general $d$-dimensional Navier-Stokes system
(NSS) is written for $d$ unknown functions
$u(t,x)=(u_1(t,x),\dots,u_d(t,x))$ of $d$ variables
$x=(x_1,\dots,x_d)$ and time $t$, and the pressure $p(t,x)$.
\begin{eqnarray}
  \frac{\partial u_j}{\partial t} + \sum_{k=1}^d u_k \frac{\partial u_j}{\partial
  x_k}&=& \nu \triangle u_j - \frac{\partial p}{\partial x_j} +
  f^{(j)} \label{eq:NS} \\
  \mbox{div}\ u&=& \sum_{j=1}^d \frac{\partial u_j}{\partial x_j}=0 \label{eq:div}
\end{eqnarray}
The functions $f^{(j)}$ are the components of the external
forcing, $\nu >0$ is the viscosity.

Therefore we have $\lambda_k = \nu |k|^2$, hence  we $p=2$ and $r=1$.

We consider (\ref{eq:NS}),(\ref{eq:div}) on the torus $\mathbb{T}^d=\left({\mathbb{R}/2\pi}\right)^d$. This enables us to use Fourier
series. We write
\begin{equation}
  u(t,x)=\sum_{k \in \mathbb{Z}^d} u_k(t)e^{\mathrm{i}(k,x)}, \qquad
  p(t,x)=\sum_{k \in \mathbb{Z}^d} p_k(t)e^{\mathrm{i}(k,x)}
\end{equation}
Observe that $u_k(t) \in \mathbb{C}^d$, i.e. they are
$d$-dimensional vectors and $p_k(t) \in \mathbb{C}$. We 
assume that $f_0=0$ and $u_0=0$.

Then (see \cite{Z} for details) (\ref{eq:div}) is reduced to 
\begin{eqnarray}
    (u_k,k)=0 \quad k \in \mathbb{Z}^d,  \label{eq:NSdivfourier}
\end{eqnarray}
and on the space of functions satisfying (\ref{eq:div}) the pressure disappears and we obtain the following infinite ladder of differential equations
for $u_k$
\begin{equation}
  \frac{d u_k}{d t}=-\mathrm{i} \sum_{k_1}(u_{k_1}|k)\sqcap_k u_{k-k_1} - \nu k^2u_k + \sqcap_k f_k,
    \label{eq:NSgal1}
\end{equation}
where $f_k$ are components of the external forcing, $\sqcap_k$
denotes the operator of orthogonal projection onto the
$(d-1)$-dimensional plane orthogonal to $k$.

Observe that the subspace defined by incompressibility condition (\ref{eq:NSdivfourier}) and reality condition (\ref{eq:reality})  is invariant under  (\ref{eq:NSgal1})
and also is invariant for  Galerkin projections of (\ref{eq:NSgal1}), where for all $k \in \mathbb{Z}^d$ all components of $u_k$ and $u_{-k}$ are both included or excluded
in the projection.  This  defines a Galerkin filtration. Alternatively, using some more specific boundary conditions  (see for example \cite{AK21,BB21}), we define a set independent modes $\{u_k\}_{k \in I}$  (maybe with respect to some other function basis consisting of $\sin(kx)$ etc)
and rewrite the convolution term in (\ref{eq:NSgal1}) using only $u_k, k \in I$.  The formulas will be a bit more complicated, but essentially still could be seen as some
convolutions.

\subsection{Preparatory remarks}
Consider equation (\ref{eq:genpde}). We assume that $N$ is a polynomial and by $D^s u$ we denote $s^{\text{th}}$ order derivative of $u$, i.e. the collection of all spatial (i.e. with respect to variable $x$) partial derivatives of $u$ of order $s$. The reason to consider polynomial and not more general functions $N$ is technical --- we need to compute the Fourier coefficients of $N\left(u,Du,\dots,D^ru\right)$. This can be achieved by taking suitable convolutions of Fourier expansions of $u$ and its spatial partial derivatives. For analytic $N$ the results are a bit more involved, i.e. we will have infinite series of convolutions. This still it is manageable but we omit it for the sake
of simplicity.

In what follows all considerations will be done assuming $u_k$ represent the Fourier expansion of $u$,
\begin{equation}
  u(t,x)=\sum_{k \in \mathbb{Z}^d} u_k(t) \exp(\mathrm{i} k \cdot x), \label{eq:u-furier}
\end{equation}
 and $u_k$ are independent variables. The estimates
and results obtained under this assumptions can be easily translated to the case when some constraints coming from the reality requirement (\ref{eq:reality}), boundary conditions
lead us to use $\sin(kx)$, $\cos(kx)$ or combinations  thereof as the basis for the expansion.  In such situation we have a natural imbedding into the previous case of the Fourier expansion with respect to
$\exp(\mathrm{i} kx)$ and all estimates are the same up to some constants.

In the sequel for $k \in \mathbb{Z}^d$ we will use norm $|k|$ which should satisfy \newline
 $|(k_1,\dots,k_d)| \geq |k_j|$ for $j=1,\dots,d$. Moreover, $u_k$ can be vectors or complex numbers. This can be easily reformulated to fit into the framework discussed in previous section, we will use a countable set of indices  $(k,\ell,a)$, where $k \in I$, $\ell \in \{1,\dots,d\}$, $a \in \{0,1\}$, so that $u_{k,\ell,a}$ is real ($a=0$) or imaginary ($a=1$) part of $\ell$-th component of vector $u_k \in \mathbb{C}^d$. The Galerkin filtration could be chosen as follows 
 \begin{equation}
 J_n=\{(k,\ell,a), |k| \leq n, \ell \in \{1,\dots,d\}, a \in \{0,1\}\}.
 \end{equation}

 After formally inserting the Fourier expansion   for $u,Du,\dots,D^r u$  in $N()$ we obtain a sum
 expressions of the following type for each monomial in $N$
\begin{equation*}
  \sum_{k_1+\dots + k_l=k} v_{k_1} \cdot  v_{k_2} \cdot
  \dots \cdot v_{k_l},
\end{equation*}
where  each of the variables $v_{k_j}$, $j=1,\dots,l$ is some Fourier coefficient of one of the components of $u$ or its partial derivatives of the order less than or equal to $r$.  This is a formal expansion, the questions of convergence and differentiability are treated in the following sections.

\subsection{Derivatives with respect to Fourier coefficients}

The goal of this section it to write  compact formulas for derivatives with respect of Fourier coefficients.

We will use the following notation to denote (partial) derivatives. For a function $f(x_1,\dots,x_n)$ and $\alpha \in \mathbb{N}^n$ we set  $D^\alpha f(x_1,\dots,x_n)=\frac{\partial^{|\alpha|}}{\partial x_1^{\alpha_1} \partial x_2^{\alpha_2}\dots \partial x_n^{\alpha_n}}$, where $|\alpha|=\sum_j \alpha_j$.
For $k \in \mathbb{Z}^n$ we define $k^\alpha= k_1^{\alpha_1} k_2^{\alpha_2} \cdot \dots \cdot k_n^{\alpha_n}$.

When considering $N(u,Du,\dots)$ it is convenient to name the arguments of $N$ by $u$, $D^\alpha u$, which we will use below. For example, for $N(u,u_x)=uu_x$ we will
have $\frac{\partial N}{\partial u}(u,u_x)=u_x$ and  $\frac{\partial N}{\partial u_x}(u,u_x)=u$. We will use this convention below.

We have
\begin{multline*}
   N(u+h,\dots,D^r(u+h))= N(u,\dots,D^ru)
    + \sum_{0\leq |\alpha| \leq r } \frac{\partial N}{\partial (D^\alpha u)}(u,\dots,D^ru)D^\alpha h \\
   + \sum_{\substack{\alpha_1,\alpha_2\\  0\leq |\alpha_1|\leq r,\\ 0 \leq |\alpha_2|\leq r}} \frac{1}{2} \frac{\partial^2 N}{\partial (D^{\alpha_1}u) \partial(D^{\alpha_2}u)}(u,\dots,D^ru) D^{\alpha_1}h D^{\alpha_2}h\\
   + O(\|h\|^3+\ldots+\|D^rh\|^3)
\end{multline*}
Therefore we  have
\begin{eqnarray}
  \frac{\partial N_k}{\partial u_j} = \sum_{|\alpha| \leq r}\left(\left(\frac{\partial N}{\partial (D^\alpha u)}(u,\dots,D^ru)\right)_{k-j}\mathrm{i}^{|\alpha|} j^\alpha\right)  \label{eq:derNk-f}
\end{eqnarray}
and
\begin{eqnarray}
  \frac{\partial^2 N_k}{\partial u_{j_1}  \partial u_{j_2}} = \left(\sum_{\substack{\alpha_1,\alpha_2\\  0\leq |\alpha_1|\leq r, \\ 0 \leq |\alpha_2|\leq r}}\left(\frac{\partial^2 N}{\partial(D^{\alpha_1}u)  \partial(D^{\alpha_2} u) }(u,\dots,D^ru)\right)_{k-j_1-j_2}\mathrm{i}^{|\alpha_1|+|\alpha_2|}j_1^{\alpha_1} j_2^{\alpha_2}\right). \label{eq:der2Nk-f}
\end{eqnarray}
Formulas (\ref{eq:derNk-f},\ref{eq:der2Nk-f}) after formally inserting Fourier expansions define us formal power series, which we denote by symbols $ \frac{\partial N_k}{\partial u_j}$, $\frac{\partial^2 N_k}{\partial u_{j_1}  \partial u_{j_2}}$ without actually claim that these functions represent partial derivatives. This fact is proven later for
arguments belong to some sets with good convergence properties.

\subsection{Estimates for sets with polynomial decay}

Throughout this section a polymial $N(u,\dots,D^r u)$ is fixed  and various estimates given below will obviously depend on this polynomial, but this will not be clearly
indicated.

In this subsection our goal is to prove the following bounds on the vector field induced in Fourier domain and its derivatives on set with polynomial decay.
\begin{theorem}
\label{thm:Dgen}
 Let $s > s_0=d+r$.
 If $|a_k| \leq C/|k|^s$, $|a_0| \leq C$, then  the formal series defining $N_k$, $\frac{\partial N_k}{\partial a_j}$, $\frac{\partial^2 N_k}{\partial a_{j_1} \partial a_{j_2}}$
 are absolutely convergent and
there exists $D=D(C,s)$, $D_1=D_1(C,s)$, $D_2=D_2(C,s)$
\begin{eqnarray}
  |N_k| &\leq& \frac{D}{|k|^{s-r}}, k \neq 0, \qquad \mbox{and} \qquad |N_0| \leq D,  \label{eq:Nk-estm} \\
  \left|\frac{\partial N_k}{\partial a_j}\right| &\leq& \frac{D_1 |j|^r}{|k-j|^{s-r}},  k \neq j, \qquad \mbox{and} \nonumber \\
     & &  \qquad \left|\frac{\partial N_k}{\partial a_j}\right|  \leq D_1|j|^r, k=j \label{eq:derNk-estm} \\
   \left|\frac{\partial^2 N_k}{\partial a_{j_1} \partial a_{j_2}}\right| &\leq& \frac{D_2 |j_1|^r |j_2|^r}{|k-j_1-j_2|^{s-r}},   k \neq j_1+j_2, \qquad \mbox{and} \nonumber \\
     & & \qquad   \left|\frac{\partial^2 N_k}{\partial a_{j_1} \partial a_{j_2}}\right|  \leq D_2 |j_1|^r |j_2|^r, k=j_1+j_2.  \label{eq:derNk2-estm}
\end{eqnarray}
\end{theorem}
Assertion (\ref{eq:Nk-estm}) has been proven as Lemma 3.1 in \cite{ZKS3}.
Before the proof of Theorem~\ref{thm:Dgen} we need to establish several short lemmas. Some of them are taken from Section 3.1 in \cite{ZKS3}, but we include
their short proofs for the sake of completeness.

\begin{lemma} \cite[Lemma 3.2]{ZKS3}
\label{lem:powineq}
  Let $\gamma > 1$. For any $a,b \geq 0$ the following inequality
  is satisfied
\begin{equation*}
  (a+b)^\gamma \leq 2^{\gamma-1}(a^\gamma + b^\gamma). 
\end{equation*}
\end{lemma}
\textbf{Proof:} This is an easy consequence of the convexity of
function $x \mapsto x^\gamma$ for $\gamma >1$. Namely
\begin{eqnarray*}
  (a+b)^\gamma = 2^\gamma \left(\frac{a + b}{2}\right)^\gamma \leq 2^\gamma \left(\frac{a^\gamma +
  b^\gamma}{2}\right)= 2^{\gamma-1} (a^\gamma + b^\gamma).
\end{eqnarray*}
\qed

The following  lemma was proved in \cite{Sa}
\begin{lemma} \cite[Lemma 3.3]{ZKS3}
\label{lem:estmQ} Assume that $\gamma > d$.  Then there exists
$S_Q(d,\gamma) \in \mathbb{R}$ such that for any $k \in
\mathbb{Z}^d \setminus \{0\}$ holds
\begin{equation*}
  \sum_{k_1 \in \mathbb{Z}^d \setminus \{0, k \}} \frac{1}{|k_1|^\gamma |k-k_1|^\gamma} \leq
    \frac{S_Q(d,\gamma)}{|k|^{\gamma}}.
\end{equation*}
\end{lemma}
\textbf{Proof:} From the triangle inequality and
Lemma~\ref{lem:powineq} we have
\begin{eqnarray*}
  \frac{|i|^\gamma}{|k-i|^\gamma |k|^\gamma} &\leq&
  \frac{\left( |k-i| + |k| \right)^\gamma}{|k-i|^\gamma |k|^\gamma} \\
 & & \leq  \frac{ 2^{\gamma -1}(|k-i|^\gamma + |k|^\gamma)}{|k-i|^\gamma
  |k|^\gamma} =
  2^{\gamma -1} \left( \frac{1}{|k|^\gamma} + \frac{1}{|k-i|^\gamma}
  \right).
\end{eqnarray*}
Hence
\begin{eqnarray*}
 \sum_{k \in {\mathbb Z}^d\setminus \{0,i\}} \frac{1}{|k|^\gamma |i -
  k|^\gamma} \leq  \sum_{k \in {\mathbb Z}^d\setminus \{0,i\}} \frac{2^{\gamma-1}}{|i|^\gamma}
  \left( \frac{1}{|k|^\gamma} +\frac{1}{ |i -
  k|^\gamma}\right) <  \\ \frac{2^{\gamma}}{|i|^\gamma} \sum_{k \in {\mathbb Z}^d \setminus \{0\}}
  \frac{1}{|k|^\gamma}.
\end{eqnarray*}
 \qed

Now we want to include  also the vectors of zero length in the sum
appearing in Lemma~\ref{lem:estmQ}.  To make
expression of some formulas less cumbersome in this subsection for
$0=\{0\}^d \in \mathbb{Z}^d$ we redefine its norm by setting
$|0|=1$.
\begin{lemma} \cite[Lemma 3.4]{ZKS3}
\label{lem:estmQ2} Assume that $\gamma > d$.  Then there exists
$C_2(d,\gamma) \in \mathbb{R}$ such that for any $k \in
\mathbb{Z}^d$ holds
\begin{equation*}
  \sum_{\substack{k_1,k_2 \in \mathbb{Z}^d\\ k_1+k_2=k}} \frac{1}{|k_1|^\gamma |k_2|^\gamma} \leq
    \frac{C_2(d,\gamma)}{|k|^{\gamma}}.
\end{equation*}
\end{lemma}
\textbf{Proof:} Consider two cases $k=0$ and $k \neq 0$.

If $k=0$, then there exists $\widetilde{C}(d,\gamma) \in
\mathbb{R}$ such that
\begin{eqnarray*}
  \sum_{\substack{k_1,k_2 \in \mathbb{Z}^d\\ k_1+k_2=k}} \frac{1}{|k_1|^\gamma |k_2|^\gamma} =
  1 + \sum_{k_1 \in  \mathbb{Z}^d\setminus \{0\}}
  \frac{1}{|k_1|^{2\gamma}} =\widetilde{C}(d,\gamma).
\end{eqnarray*}
If $k \neq 0$, then from Lemma~\ref{lem:estmQ} it follows that
\begin{eqnarray*}
 \sum_{\substack{k_1,k_2 \in \mathbb{Z}^d\\ k_1+k_2=k}} \frac{1}{|k_1|^\gamma
 |k_2|^\gamma}= \frac{2}{|k|^\gamma} +
 \sum_{\substack{k_1,k_2 \in \mathbb{Z}^d \setminus \{0\}\\ k_1+k_2=k}} \frac{1}{|k_1|^\gamma |k_2|^\gamma}
    \leq
    \frac{S_Q(d,\gamma) +   2}{|k|^\gamma}.
\end{eqnarray*}
Hence the assertion holds for
$C_2(d,\gamma)=\max(\widetilde{C}(d,\gamma),S_Q(d,\gamma) + 2 )$.
\qed

\begin{lemma} \cite[Lemma 3.5]{ZKS3}
\label{lem:estmQn} Assume $\gamma > d$. For any $n \in
\mathbb{Z}_+$, $n > 1$  there exists $C_n(d,\gamma) \in
\mathbb{R}$ such that for any $k \in \mathbb{Z}^d$ holds
\begin{equation*}
  \sum_{k_1,k_2,\dots,k_n \in \mathbb{Z}^d, \sum_{i=1}^n k_i=k}
  \frac{1}{|k_1|^\gamma |k_2|^\gamma \cdot \dots \cdot |k_n|^\gamma} \leq
    \frac{C_n(d,\gamma)}{|k|^{\gamma}}.
\end{equation*}
\end{lemma}
\textbf{Proof:} By induction. Case $n=2$ is contained in
Lemma~\ref{lem:estmQ2}. Assume now that the assertion holds for
$n$. We have
\begin{eqnarray*}
  \sum_{k_1,k_2,\dots,k_{n+1} \in \mathbb{Z}^d, \sum_{i=1}^{n+1} k_i=k}
  \frac{1}{|k_1|^\gamma |k_2|^\gamma \cdot \dots \cdot |k_{n+1}|^\gamma}
   = \\ \sum_{k_{n+1} \in \mathbb{Z}^d} \left( \frac{1}{|k_{n+1}|^\gamma}
    \sum_{k_1,k_2,\dots,k_n \in \mathbb{Z}^d, \sum_{i=1}^n k_i=k-k_{n+1}}
  \frac{1}{|k_1|^\gamma |k_2|^\gamma \cdot \dots \cdot |k_n|^\gamma}
  \right) \leq \\
 \sum_{k_{n+1} \in \mathbb{Z}^d}  \frac{1}{|k_{n+1}|^\gamma} \cdot
 \frac{C_n(d,\gamma)}{|k-k_{n+1}|^\gamma} \leq \frac{C_2(d,\gamma)
 C_n(d,\gamma)}{|k|^\gamma}.
\end{eqnarray*}
\qed

\noindent \textbf{Proof of Theorem~\ref{thm:Dgen}:} For the proof it
is enough to assume that $N$ is a monomial. After formally
inserting the Fourier expansion  for $u,Du,\dots,D^r u$ we obtain
the expression of the following type
\begin{equation}
  N_k(u)=\sum_{k_1+\dots + k_l=k} v_{k_1} \cdot  v_{k_2} \cdot
  \dots \cdot v_{k_l}, \label{eq:Nkmon}
\end{equation}
where  each of the variables $v_{k_i}$, $i=1,\dots,l$ is some
Fourier coefficient of one the components of $u$ or its partial
derivatives of the order less than or equal to $r$.

Observe that for the Fourier coefficients of partial derivatives
up to order $r$ we have the following estimates
\begin{equation}
  \left| \frac{\partial^{\beta_1 + \dots + \beta_l} u}{\partial x_1^{\beta_1}\dots \partial x_d^{\beta_l}
  }\right| \leq \frac{C}{|k|^{s - (\beta_1 + \dots + \beta_l)}}
  \leq \frac{C}{|k|^{s - r}}.  \label{eq:betaestm}
\end{equation}
From conditions (\ref{eq:Nkmon}) and (\ref{eq:betaestm}), and
Lemma~\ref{lem:estmQn} we obtain
\begin{equation*}
 |N_k(u)| \leq \sum_{k_1+\dots + k_n=k} \frac{C^n}{|k_1|^{s-r} \cdot \dots \cdot |k_n|^{s-r}}
   \leq \frac{C^n C_n(d,s-r)}{|k|^{s-r}}
\end{equation*}
This establishes (\ref{eq:Nk-estm}).

To prove (\ref{eq:derNk-estm}) we use (\ref{eq:derNk-f}). From previous reasoning applied to polynomials $\frac{\partial N}{\partial (D^\alpha u)}(u,\dots,D^ru)$ we immediately
obtain
\begin{eqnarray*}
 \left|\left(\frac{\partial N}{\partial (D^\alpha u)}(u,\dots,D^ru)\right)_{k-j}\right| \leq \frac{D_1}{|k-j|^{r-s}},
\end{eqnarray*}
which together with bound $|\alpha| \leq r$ gives (\ref{eq:derNk-estm}).

The proof of (\ref{eq:derNk2-estm}) is analogous. \qed

Later we will use the following lemma.
\begin{lemma}
\label{lem:sumForDer}
Assume that $s > 2r + d$, then there exists $C(d,s,r)$ such that
\begin{equation*}
   \sum_{k \in \mathbb{Z}^d, k \neq i}\frac{|k|^r}{|i-k|^{s-r}} \leq  C(d,s,r)(1+|i|^r).
\end{equation*}
\end{lemma}
\textbf{Proof:}
\begin{multline*}
    \sum_{k \in \mathbb{Z}^d, k \neq i}\frac{|k|^r}{|i-k|^{s-r}} \leq \sum_{k \in \mathbb{Z}^d, k \neq i}\frac{(|i-k| + |i|)^r}{|i-k|^{s-r}} \\
\leq \sum_{k \in \mathbb{Z}^d, k \neq i}\frac{2^{r-1}|i-k|^r + 2^{r-1}|i|^r}{|i-k|^{s-r}}\\
= 2^{r-1} \left(\sum_{k \in \mathbb{Z}^d, k \neq i}\frac{1}{|i-k|^{s-2r}} + \sum_{k \in \mathbb{Z}^d, k \neq i}\frac{|i|^r}{|i-k|^{s-r}}\right) \leq
C(d,s,r)(1 + |i|^r).
\end{multline*}
Observe that the two infinite sums are convergent under assumption $s> 2r +d$.
\qed

\subsection{Conditions S1 and S2}
\label{subsec:ch-concond}
We will be interested in the following candidates for space $H$:
\begin{itemize}
\item $c_0$, $\|x\|=\sup_{k\in\mathbb Z^d} |x_k|$,
\item $l_p$, $\|x\|=\left(\sum_{k\in\mathbb Z^d}  |x_k|^p\right)^{1/p}$ with $p \geq 1$,
\item $\mathcal{W}^{m,p}$, $\|x\|=\left( \sum_{i=0}^m  \sum_{k\in\mathbb Z^d}  |k|^{ip} |x_k|^p \right)^{1/p}$ with $p \geq 1$ and $m \geq 0$.
\end{itemize}
It is immediate that the space listed above are \gss spaces.

We consider sets of the form
\begin{eqnarray*}
  W_{P}(C,s) &=&\left\{ |x_k| \leq \frac{C}{|k|^s}, k \neq 0,  |x_0| \leq C \right\}, s >0, \\
  W_{exp}(q,S)&=&\left\{ |x_k| \leq \frac{S}{q^{|k|}} \right\}, \quad q>1.
\end{eqnarray*}
Clearly they satisfy condition \textbf{S1}.

\begin{remark}
	Observe, that if $W_{exp}(q,S)\subset W_P(C,s)\subset H$ then it is a closed subset and thus it inherits {\bf S2}, {\bf C1} and {\bf C2} from $W_P(C,s)$. Thus, in what follows we will focus on sets with polynomial decay, only.
\end{remark}

We will show that condition \textbf{S1} on sets $W_{exp}(q,S)$, $W_P(C,s)$ when considered in suitable spaces.
\begin{theorem}
\label{thm:polexp-comp}
The following statements hold true.
\begin{itemize}
 \item $W_P(C,s)$ is a compact subset of  $c_0$ for any $d$ and $s>1$.
 \item $W_P(C,s)$ is a compact subset of $l_p$ and $\mathcal{W}^{m,p}$ provided $(s-m)p >d$ (with $m=0$ for space  $l_p$).
 \item  $W_{exp}(q,S)$ is a compact subset of  $c_0$, $l_p$, $\mathcal{W}^{m,p}$ for any $d$, $m$, $p$ and $q>1$.
\end{itemize}
\end{theorem}
\textbf{Proof:}
We  use Lemma~\ref{lem:compt-gss} as a criterion for compactness. It is immediate that the sets $ W_{P}(C,s)$ and $ W_{exp}(q,S)$ are closed in all spaces considered (if contained in them).

The case of $c_0$ space is obvious. In other cases it is enough that to observe that
\begin{equation*}
  \sum_{|k| \geq n} |k|^a |x_k|^b \sim \int_{n}^\infty r^{d-1} r^a x(r)^{b}dr,
\end{equation*}
where $x(r)=\frac{C}{r^s}$ in case od $W_P$ and $x(r)=S q^{-r}$ in the case od $W_{exp}$. From this we obtain our assertion for $W_{exp}$ for all $d$, $p$ and $m$.

In the case of $W_P$ the largest exponent under integral must be less than $-1$ to obtain the convergence, which gives
\begin{equation*}
  d-1 + mp - sp < -1.
\end{equation*}
\qed

\subsection{Conditions \textbf{C1,C2}}

\subsubsection{Condition C1 for set with polynomial decay}

\begin{theorem}
\label{thm:existsCi} Consider (\ref{eq:fugenpde}). Assume that
conditions (\ref{eq:lambdak}) and
(\ref{eq:p>r}) hold.

Let $W=W_{P}(C,s)$. Let $H$ be one of space listed at the beginning of Section~\ref{subsec:ch-concond}.

Then
\begin{itemize}
  \item if $s > p+d$, then $W \subset c_0$ and $F:W \to c_0$ is continuous
  \item if $(s-p - \ell)q > d$, then   $W \subset \mathcal{W}^{\ell,q}$   and
  $F:W \to \mathcal{W}^{\ell,q}$ is continuous (with $\ell=0$ for space  $l_q$)
\end{itemize}

\end{theorem}
\textbf{Proof:}

The first question is whether  $W \subset \dom{F}$.
Consider $u \in W$. From Theorem~\ref{thm:Dgen} and condition (\ref{eq:lambdak}) it follows
that $F_k(u)$ is defined and for $|k|>K$ holds
\begin{equation}
  |F_k(u)| \leq L^* C |k|^{p-s} + D |k|^{r-s} \leq
  \frac{D_2}{|k|^{s-p}} \label{eq:fkbd}.
\end{equation}
for some constants $D$ and $D_2$.
From Theorem~\ref{thm:polexp-comp} and our assumptions it follows $W$ is compact in various spaces listed in the assertion and $F(W)$ is also in the listed space
and is contained in some compact space (due to Lemma~\ref{lem:compt-gss} and decay estimates (\ref{eq:fkbd})). Moreover for all  $u \in W$
\begin{equation*}
 \lim_{n \to \infty} P_n F(u)=F(u).
\end{equation*}
uniformly on $W$.

Hence to prove the continuity of $F:W \to H$ it is enough to prove that $F_k:W \to H_k$ is continuous.

Let us fix index $k$ and assume   $u^n,u^* \in W$, for $n
\in \mathbb{N}$ and $u^n \to u^*$ for $n \to \infty$. We have
(compare the proof of Theorem~\ref{thm:Dgen})
\begin{equation*}
  F_k(u)=\lambda_k u_k + N_k(u)=\lambda_k u_k + \sum_{i \in J} N_{k,i}(u),
\end{equation*}
where $J$ is some set of multindices and for each $i \in J$,
$N_{k,i}$ is monomial depending on the finite number of $u_{l}$,
i.e.
\begin{displaymath}
  N_{k,i}= a u_{k_1}\cdot u_{k_2}  \cdot \dots u_{k_l}, \quad
  \mbox{for some $a \in \mathbb{C}$ and $k_1+\dots+k_l=k$ }
\end{displaymath}

The term $\lambda_k u_k$ is continuous, hence it is enough to
consider $N_k$, only. Let us fix $\epsilon > 0$. From
Theorem~\ref{thm:Dgen} it follows that there exists a finite set $S
\subset J$, such that
\begin{equation}
  \sum_{i \in J \setminus S} |N_{k,i}(u)| < \epsilon/3, \qquad
    \mbox{ for all $u \in W$}. \label{eq:Ntailestm}
\end{equation}
There exists $L$, such that for all $i \in S$ monomials
$N_{k,i}(u)$ depend in fact on the variables $u_{l}$ for $|l| \leq
L$, hence $\sum_{i \in S} N_{k,i}(u)$ is continuous on $W$. Therefore there exists $n_0$, such that
\begin{equation}
  \left| \sum_{i \in S} N_{k,i}(u^n) -  \sum_{i \in S} N_{k,i}(u^*)
  \right|< \epsilon/3.  \label{eq:Nmainestm}
\end{equation}
From (\ref{eq:Nmainestm}) and (\ref{eq:Ntailestm}) we obtain for
$n > n_0$
\begin{eqnarray*}
  |N_k(u^n) - N_k(u^*)| \leq   \left| \sum_{i \in S} N_{k,i}(u^n) -  \sum_{i \in S} N_{k,i}(u^*)
  \right| +  \\ \sum_{i \in J \setminus S} |N_{k,i}(u^n)| +
     \sum_{i \in J \setminus S} |N_{k,i}(u^*)| < \epsilon.
\end{eqnarray*}
Hence $\lim_{n \to \infty} N_k(u^n)=N_k(u^*)$.
 \qed

\subsubsection{Condition C2 for sets with polynomial decay}
\label{sssec:condDpolDec}
We will treat only the spaces $l_2$, $c_0$ and $l_1$ because for those spaces we have nice and relatively compact formulas for logarithmic norms.

Namely, we have for $A \in \mathbb{R}^{n \times n}$
\begin{eqnarray*}
  \mu_2 (A) = \max (\lambda \in Sp(A+A^t)/2), \\
  \mu_1 (A) = \max_{i=1,\dots,n} \left(a_{ii} + \sum_{k,k \neq i} |a_{ik}| \right), \\
  \mu_\infty(A) = \max_{k=1,\dots,n} \left(a_{kk} + \sum_{i,k \neq i} |a_{ik}| \right).
\end{eqnarray*}
\begin{lemma}
The same assumptions and notation as in Lemma~\ref{thm:existsCi}, but we restrict $H$ to one of the following spaces $c_0$, $l_2$, $l_1$. Assume additionally that $s>2r+d$. Then $F$ satisfies
condition {\bf C2} on $W$.
\end{lemma}
\textbf{Proof:}
First we deal with $H=c_0$. In this case the logarithmic norm is given by $\mu_\infty$.
From Theorem~\ref{thm:Dgen} and Lemma~\ref{lem:sumForDer} we have on $W$
\begin{eqnarray*}
S(row)_k:=\sum_{k,k \neq i} \left|\frac{\partial F^n_{i}}{\partial a_k} \right| \leq  \sum_{k \in \mathbb{Z}^d, k \neq i}\frac{D_1|k|^r}{|i-k|^{s-r}} \leq \tilde{C}_1(1 + |i|^r).
\end{eqnarray*}
Since  by (\ref{eq:lambdak}) $a_{kk}   \leq -L_* |k|^p$, hence since $p>r$ (by (\ref{eq:p>r})) we obtain \\
 $S(row)_k + a_{kk} \to -\infty$ for $|k| \to \infty$.
Therefore there exists $l \in \mathbb{R}$, such that
$\mu_{\infty}(Df^n(P_nz))< l$ for all $n$ and $z \in W$.

Now we assume  $H=l_1$ and the logarithmic norm is $\mu_1$. We have
\begin{eqnarray*}
S(col)_i:=\sum_{i,k \neq i} \left|\frac{\partial F^n_{i}}{\partial a_k} \right| \leq  \sum_{i \in \mathbb{Z}^d, k \neq i}\frac{D_1|k|^r}{|i-k|^{s-r}} \leq  |k|^r \tilde{C}_2.
\end{eqnarray*}
We conclude exactly as in the case of $\mu_\infty$.

When $H=l_2$, then we use Gershogorin Theorem \cite{G} to estimate the spectrum $\frac{1}{2}\left(DF^n + \left(DF^n\right)^T\right)$. Observe the radius of $k$-th Gershgorin circle
will be bounded from above as follows
\begin{equation*}
  R_i \leq (S(col)_i + S(row)_i)/2 \leq \tilde{C}|i|^r.
\end{equation*}
We conclude as in the previous cases.
\qed

%% file: varconvcond.tex
\subsection{Convergence conditions for variational system}

Let us fix $W=W_P(C,s)$ with \textbf{S1}, \textbf{S2} and such that the vector field $F=L+N$ satisfies \textbf{C1} and \textbf{C2} on $W$. On variational variables we also take set
of the form $W_V=W_P(C_V,s)$. We will show that $F$ also satisfies \textbf{VC} on $W\times W_V$.

From Theorem~\ref{thm:Dgen} we have for $z \in W$
\begin{eqnarray}
   \left|\frac{\partial N_k}{\partial a_j}(z)\right| &\leq& \frac{D_1 |j|^r}{|k-j|^{s-r}},  k \neq j, \qquad \mbox{and} \nonumber \\
     & &  \qquad \left|\frac{\partial N_k}{\partial a_j}(z)\right|  \leq D_1|j|^r, k=j \label{eq:derNk-e-variso}
\end{eqnarray}

\subsubsection{Condition \textbf{VC1} for sets with polynomial decay}
\begin{lemma}
The same assumptions and notation as in Theorem~\ref{thm:existsCi}.  Assume that $s>r+d$ and if $H=l_q$ or $H=\mathcal{W}^{m,q}$, then $(s-r - m)q>d$ (with $m=0$ for space $l_q$).

Then $F$ satisfies condition {\bf VC1} on $W \times W_V$.
\end{lemma}
\textbf{Proof:}
From (\ref{eq:derNk-e-variso}) and  Lemma~\ref{lem:estmQ} it follows that for $(z,\mathcal{C}) \in W\times W_V$ holds
\begin{eqnarray*}
  |\widetilde{DF}_i(u) \mathcal{C}| \leq \sum_j \frac{D_1|j|^r}{|i-j|^{s-r}} \frac{C_V}{|j|^s}=D_1 C_V\sum_j \frac{1}{|i-j|^{s-r}|j|^{s-r}}\\
  =D_1C_VS_Q(d,s-r)\frac{1}{|i|^{s-r}}
\end{eqnarray*}
From considerations as in the proof of Theorem~\ref{thm:polexp-comp} it follows that  the image of $\widetilde{DF}(z) \mathcal{C}$ will be contained in compact set $W_P(D_1C_VS_Q(d,s-r),s-r)$ in the following cases
\begin{itemize}
\item $H=c_0$ if $s-r >1$
\item $H=l_q$ or $H=\mathcal{W}^{m,q}$ provided $(s-r - m)q>d$ (with $m=0$ for space $l_q$).
\end{itemize}
Therefore for continuity of map $F_V(u,\mathcal{C})$ it is enough that its Galerkin projections are continuous (which is obvious - see the proof of continuity of $N$
in Theorem~\ref{thm:existsCi}), because then in view of the above
estimates they converge to $F_V$ on $W \times W_V$.
\qed

\subsubsection{Condition \textbf{VC2} for sets with polynomial decay}

\begin{lemma}
The same assumptions and notation as in Theorem~\ref{thm:existsCi}, but we restrict $H$ to one of the following spaces $c_0$, $l_2$, $l_1$. Assume additionally that $s>2r+d$. Then $F$ satisfies
condition {\bf VC2} on $W \times W_V$.
\end{lemma}
\textbf{Proof:}
Observe that
\begin{equation*}
  DP_nF_{V|H_n \times H_n}(x,\mathcal{C})= \left[\begin{array}{cc}
                                         DF^n(x)  & 0 \\
                                         D^2F^n(x)\mathcal{C} & DF^n(x)
                                        \end{array} \right]
\end{equation*}
To be more precise we have for  $|i|,|j| \leq n$
\begin{eqnarray*}
   \left(DP_nF_{V|H_n \times H_n}(x,\mathcal{C})\right)_{x_i,x_j}&=&\frac{\partial F^n_i}{\partial x_j},\\
    \left(DP_nF_{V|H_n \times H_n}(x,\mathcal{C})\right)_{x_i,\mathcal{C}_j}&=&0,  \\
     \left(DP_nF_{V|H_n \times H_n}(x,\mathcal{C})\right)_{\mathcal{C}_i,x_j}&=&\sum_{|k|\leq n}\frac{\partial^2 F^n_i}{\partial x_j\partial x_k}\mathcal{C}_k,  \\
     \left(DP_nF_{V|H_n \times H_n}(x,\mathcal{C})\right)_{\mathcal{C}_i,\mathcal{C}_j}&=&\frac{\partial F^n_i}{\partial x_j}.
\end{eqnarray*}
Formulas for logarithmic norms in $c_0$, $l_1$, $l_2$ has been given in section~\ref{sssec:condDpolDec}.  From these formulas it is clear that to be establish condition
\textbf{VC2} it is enough to show that on $W \times W_V$ the following holds for some constants $E_1$, $E_2$
\begin{eqnarray}
  \sum_{|j| \leq n}\left|   \left(DP_nF_{V|H_n \times H_n}(x,\mathcal{C})\right)_{\mathcal{C}_i,x_j} \right| \leq E_1 |i|^r,  \label{eq:var-row-sum} \\
  \sum_{|i| \leq n}\left|   \left(DP_nF_{V|H_n \times H_n}(x,\mathcal{C})\right)_{\mathcal{C}_i,x_j} \right| \leq E_2|j|^r.  \label{eq:var-col-sum}
\end{eqnarray}

We start with (\ref{eq:var-row-sum}). From Theorem~\ref{thm:Dgen} and Lemmas~\ref{lem:estmQ}  and \ref{lem:sumForDer}
\begin{eqnarray*}
  \sum_{|j| \leq n}\left|   \left(DP_nF_{V|H_n \times H_n}(x,\mathcal{C})\right)_{\mathcal{C}_i,x_j} \right| \leq \sum_{|j|\leq n}  \sum_{|k|\leq n}\left|\frac{\partial^2 F^n_i}{\partial x_j\partial x_k}\mathcal{C}_k \right| \\
  = \sum_{|j|\leq n}  \sum_{|k|\leq n}\left|\frac{\partial^2 N_i}{\partial x_j\partial x_k}\mathcal{C}_k \right| \leq
  \sum_{|j|\leq n}  \sum_{|k|\leq n}  \frac{D_2 |j|^r |k|^r}{|i-j-k|^{s-r}} \frac{C_V}{|k|^s} \\
    \leq  D_2 C_V \sum_{j}  \sum_{k}  \frac{|j|^r }{|i-j-k|^{s-r} |k|^{s-r}} \\
    \leq D_2 C_V S_Q(d,s-r)  \sum_{j} \frac{|j|^r}{|i-j|^{s-r}} \leq D_2C_V S_Q(d,s-r)  C(d,s,r)(1+|i|^r).
\end{eqnarray*}

Now we look establish (\ref{eq:var-col-sum}). From Theorem~\ref{thm:Dgen} and Lemma~\ref{lem:estmQ}   we obtain
\begin{eqnarray*}
  \sum_{|i| \leq n}\left|   \left(DP_nF_{V|H_n \times H_n}(x,\mathcal{C})\right)_{\mathcal{C}_i,x_j} \right| \leq \sum_{|i|\leq n} \sum_{|k|\leq n}\left|\frac{\partial^2 F^n_i}{\partial x_j\partial x_k}\mathcal{C}_k\right| \\
  =\sum_{|i|\leq n} \sum_{|k|\leq n}\left|\frac{\partial^2 N_i}{\partial x_j\partial x_k} \right| \cdot  |\mathcal{C}_k|  \leq \sum_{|i|\leq n} \sum_{|k|\leq n} \frac{D_2 |j|^r |k|^r}{|i-j-k|^{s-r}} \frac{C_V}{|k|^s} \\
  \leq D_2 C_V \sum_{i} \sum_{k} \frac{|j|^r }{|i-j-k|^{s-r} |k|^{s-r}} \\
  \leq  D_2 C_V S_Q(d,s-r) |j|^r \sum_{i}  \frac{1}{|i-j|^{s-r}}\\
  =D_2 C_V S_Q(d,s-r)\left( \sum_{i \in \mathbb{Z}^d}  \frac{1}{|i|^{s-r}} \right) |j|^r
\end{eqnarray*}
\qed

\subsection{Isolation property for the main system }
\label{subsec:iso-main-var}

\begin{lemma}
	\label{lem:iso-mainVar}
	The vector field (\ref{eq:fugenpde}) is isolating on the family of sets
	$$\mathcal W = \left\{W_P(C,s) : C>0,\ s>d+r\right\}.$$
\end{lemma}
\textbf{Proof:}
Fix $C>0$ and $s>d+r$ and $W=W_P(C,s)$. From Theorem~\ref{thm:Dgen} it follows that there exists
$D=D(C,s)$, such that
\begin{equation*}
	|N_k(u)| < \frac{D}{|k|^{s - r}}, \quad \mbox{for all $u \in W$.}
\end{equation*}
Let $L_*>0, p>r, K$ be constants as in (\ref{eq:lambdak}).
Then for $u \in W $ and all $|k|>K$ such that $|u_{k}|=\frac{C}{|k|^{s}}$ there holds
\begin{multline*}
	\frac{1}{2}\frac{d}{dt}(u_{k}|u_{k}) < -L_* |k|^p |u_{k}|^2 + |u_{k}| |N_{k}(u)| \leq  \\
	\left(- L_* C |k|^{p-s} + D |k|^{r -s}\right) |u_{k}| a
	= \left(- L_* C |k|^{p-r} + D \right) |k|^{r-s}|u_{k}|.
\end{multline*}
Since $p>r$ we conclude that $\frac{d}{dt}\|u_k\|^2<0$ for $|k|>K_0$ sufficiently large. The same inequalities hold true for all Galerkin projections of (\ref{eq:fugenpde}) for $n>K_0$.
\qed

\subsection{Isolation property for variational equation}

We want to show the isolation property for the variational part of system (\ref{eq:sysVar1})-(\ref{eq:sysVar2}).

\begin{lemma}
\label{lem:iso-VarVar}
The vector field $F_V$ induced by (\ref{eq:fugenpde}) is isolating on the family of sets
$$
\mathcal W =\left\{W_P(C,s)\times W_P(C_V,s) : C>0,\ C_V>0,\ s>d+r\right\}.
$$
Moreover, the constant $K_0$ does not depend on $C_V$.
\end{lemma}
\textbf{Proof:} Fix $W\times W_V\in\mathcal W$. In view of Lemma~\ref{lem:iso-mainVar} we just have to prove
\begin{equation}
	\frac{1}{2}\frac{d}{dt}(\mathcal{C}_k,\mathcal{C}_k) < 0, \quad \mbox{if $|\mathcal{C}_k| = \frac{C_V}{|k|^s}$}. \label{eq:iso-VarVar}
\end{equation}

Let $L_*>0, p>r, K$ be constants as in (\ref{eq:lambdak}). From (\ref{eq:derNk-e-variso}) for any $(z,V) \in W\times W_V$ there holds (we use Lemma~\ref{lem:estmQ} and the assumption $s-r>d$)
\begin{eqnarray*}
  \left|\sum_j \frac{\partial N_k}{\partial a_j}(z)V_j \right| \leq  D_1 C_V \sum_j \frac{|j|^r}{|k-j|^{s-r} |j|^{s}}\\
  =D_1 C_V \sum_j \frac{1}{|k-j|^{s-r} |j|^{s-r}}=
  D_1 C_V S_Q(d,s-r) \frac{1}{|k|^{s-r}}.
\end{eqnarray*}
From this for $|k|>K$ and $z\in W$ and $\mathcal C_k\in\mathrm{bd}(\pi_k W_V)$ we get
\begin{multline*}
	\frac{1}{2}\frac{d}{dt}(\mathcal{C}_k,\mathcal{C}_k)  < \lambda_k\frac{C_V}{|k|^s} + D_1 C_V S_Q(d,s-r) |k|^{r-s}\\
	< -L_* C_V|k|^{p-s} + D_1 C_V S_Q(d,s-r)|k|^{r-s} \\
	= \left(-L_*|k|^{p-r} + D_1 S_Q(d,s-r) \right)C_V|k|^{r-s}.
\end{multline*}
Hence, for
\begin{equation*}
	|k| > K_0:=\left(\frac{D_1 S_Q(d,s-r)}{L_*} \right)^{\frac{1}{p-r}}
\end{equation*}
we obtain (\ref{eq:iso-VarVar}).
\qed

\begin{remark}
For our further investigations it is crucial that the constant $K_0$ in Lemma~\ref{lem:iso-VarVar} does not depend on $C_V$. This is a natural consequence of linearity of equations with respect to variational variable.
\end{remark}